\documentclass[AMA,STIX1COL]{WileyNJD-v2}

\articletype{Research Article}%

\received{...}
\revised{...}
\accepted{...}

\usepackage{xcolor}
\usepackage{amsfonts}
\usepackage{graphicx}
\usepackage{epstopdf}
\usepackage{overpic}
\usepackage{amsmath,esint}
\usepackage{upquote}
\usepackage{extarrows}
\usepackage{overpic}
\usepackage{amsopn}

\renewcommand{\v}{\bs{v}}
\newcommand{\bs}{\boldsymbol}
\renewcommand{\d}{\,\textup{d}}
\usepackage{tikz}
\usetikzlibrary{shapes,arrows}
\usepackage{hyperref}
\usepackage[capitalize]{cleveref}
\hypersetup{hidelinks}

\ifpdf
  \DeclareGraphicsExtensions{.eps,.pdf,.png,.jpg}
\else
  \DeclareGraphicsExtensions{.eps}
\fi

\begin{document}

\title{An optimal complexity spectral method for Navier--Stokes simulations in the ball}

\author[1]{Nicolas Boull\'e}

\author[2]{Jonasz S\l omka}

\author[3]{Alex Townsend}

\authormark{Boull\'e et al.}

\address[1]{\orgdiv{Mathematical Institute}, \orgname{University of Oxford}, \orgaddress{Oxford, \country{UK}}}

\address[2]{\orgdiv{Institute for Environmental Engineering, Department of Civil, Environmental and Geomatic Engineering}, \orgname{ETH Z\"urich}, \orgaddress{Z\"urich, \country{Switzerland}}}

\address[3]{\orgdiv{Department of Mathematics}, \orgname{Cornell University}, \orgaddress{Ithaca, \state{NY}, \country{USA}}}

\corres{Nicolas Boull\'e, Mathematical Institute, University of Oxford, Oxford, UK.\\ \email{boulle@maths.ox.ac.uk}}

\abstract[Summary]{
We develop a spectral method for solving the incompressible generalized Navier--Stokes equations in the ball with no-flux and prescribed slip boundary conditions. The algorithm achieves an optimal complexity per time step of $\mathcal{O}(N\log^2(N))$, where $N$ is the number of spatial degrees of freedom. The method relies on the poloidal-toroidal decomposition of solenoidal vector fields, the double Fourier sphere method, the Fourier and ultraspherical spectral method, and the spherical harmonics transform to decouple the Navier--Stokes equations and achieve the desired complexity and spectral accuracy.
}


\maketitle

\section{Introduction}\label{sec_intro}

Complex fluids, such as active fluids in biology and quantum fluids in physics, are often modeled by Navier--Stokes-like equations. Computational techniques for simulating these fluids in realistic configurations found in laboratories, such as cylinders, balls, and ellipsoids, are critical for experimental design. In this paper, we develop a new optimal-complexity spectral method for their simulation.  In particular, we consider the generalized Navier--Stokes (NS) equations, which have been used to model active fluids~\cite{slomka2017spontaneous,linkmann2019phase}:
\begin{subequations} \label{eq_active_fluids}
\begin{align}
\nabla\cdot \v &= 0,\\
\partial_t \v + \v\cdot \nabla \v &= -\nabla p + \nabla \cdot \bs{\sigma},
\end{align}
\end{subequations}
defined on the unit ball $\Omega=\overline{B(0,1)}$ with velocity boundary conditions tangent to the surface of the sphere, i.e.~$\v|_{\partial \Omega}\cdot\hat{\mathbf{r}} = 0$, where $\hat{\mathbf{r}}$ is the unit radial direction. Here, $\v(t,x)$ is the fluid velocity and $p(t,x)$ the local pressure, where $t\geq 0$ and $x\in \overline{B(0,1)}$. In this work, we consider the stress tensor $\bs{\sigma}(t,x)$ to be defined as
\begin{equation}
\bs{\sigma} = (\Gamma_0 - \Gamma_2\nabla^2 + \Gamma_4\nabla^4)[\nabla \v + (\nabla \v)^\top],
\end{equation}
with higher order derivatives $\nabla^{2n} = (\nabla^2)^n$ for $n\geq 2$ stimulating non-Newtonian effects. When $\Gamma_2=\Gamma_4=0$, \cref{eq_active_fluids} reduces to the incompressible NS equations with Reynolds number $Re = 1/\Gamma_0$.

Standard approaches for numerically solving the incompressible NS equations use primitive variables to solve for the fluid's velocity and pressure directly. One of the main difficulties with this formulation is ensuring that the computed velocity vector remains exactly, or even approximately, divergence-free. Typically, projection-based schemes are used to enforce this condition~\cite{chorin1967numerical,chorin1968numerical}. Still, these can lead to difficulties with boundary conditions and reduce temporal accuracy~\cite{guermond2006overview}. A novel approach is taken by the Dedalus software project, where there full-ball spectral code implicitly imposes incompressibility of the fluid in the coupled formulation~\cite{burns2020dedalus,lecoanet2018tensor} and ensures regularity of the solution using one-side Jacobi bases~\cite{vasil2018tensor}.  Another approach that appears in numerous spectral solvers~\cite{boronski2007poloidal,boronski2007poloidal2,marti2014full,matsui2016performance},  and the approach that we use in this paper, is to reformulate the NS equations using a poloidal-toroidal (PT) decomposition. This technique reduces the fully coupled NS system for the velocity and pressure to two coupled equations for the poloidal and toroidal scalar fields and has numerous computational benefits: (i) There can be a significant computational saving as only scalar fields need to be computed; (ii) Vector-valued PDEs are reduced to scalar-valued PDEs, allowing for fast PDE solvers and a parallel implementation; and (iii) The poloidal-toroidal decomposition~\cite{backus1986poloidal,chandrasekhar1961hydrodynamic,schmitt1992decomposition} explicitly enforces incompressibility of the fluid velocity vector without the need for projection-based methods. One of our novelties, is to impose the boundary conditions for the velocity using equivalent integral conditions~\cite{quartapelle1981projection} on the poloidal and toroidal scalars of the vorticity, which allows us to completely decouple the equations and achieve optimal complexity of the spectral solver at each time-step with respect to the spatial discretization.

In this paper, we describe an algorithm for solving the generalized NS equations, which achieves both spectral accuracy and optimal complexity per time step $\mathcal{O}(N\log^2(N))$ with respect to the number of degrees of freedom $N$. We refer to optimal complexity as an asymptotically linear complexity up to polylogarithmic factors of $N$. In addition to the poloidal-toroidal decomposition mentioned above, the numerical method relies on three ideas for discretizing the equations: (1) The extension of the double Fourier sphere (DFS) method to the ball~\cite{boulle2020computing}, which allows for fast computations with functions and vector fields such as vector calculus operations. (2) A spherical harmonics basis for solving the partial differential equations resulting from the PT decomposition of the NS equations. (3) The ultraspherical spectral method~\cite{olver2013fast} for performing optimal complexity mathematical operations (such as differentiation and solving PDEs) with respect to the radial variable, which relies on the relations between Chebyshev and ultraspherical polynomials. There are several significant benefits to using global spectral methods for~\cref{eq_active_fluids}: (a) Superior approximations of high-order derivatives that appear in~\cref{eq_active_fluids}, (b) better spatial resolution of advection-dominated fluids, (c) more accurate long-time simulations (by avoiding artificial diffusion), and (d) robust bulk-boundary fluid interactions. 

This paper is organized as follows. In \cref{sec_PT_NS}, we use the poloidal-toroidal decomposition to decouple the NS equations and enforce boundary conditions. Then, in \cref{sec_numer_method}, we describe the spectral discretization and associated numerical solvers for solving the poloidal and toroidal equations. Finally, we provide numerical examples in \cref{sec_numer_examples} and conclude in \cref{sec_concl}.

\section{Poloidal-toroidal formulation of the Navier--Stokes equations}
\label{sec_PT_NS}

\subsection{Vorticity formulation}
\label{sec_vor_form}

We first express the NS equations in the vorticity-vector potential form. According to the Helmholtz-Hodge decomposition~\cite{bhatia2013helmholtz}, any divergence-free vector field $\v$ in the ball can be decomposed as
\[\v=\nabla\times \bs{\psi}+\nabla \phi,\]
where the divergence-free vector potential $\bs{\psi}$ is normal to the boundary, $\bs{\psi}\times\hat{\mathbf{r}}|_{r=1}=0$, and $\phi$ is a harmonic function. Note that no harmonic vector fields are present in the above formula since our domain is simply connected. The gradient $\nabla\phi$ is uniquely determined by the Neumann boundary conditions corresponding to the flux across the boundary $\v\cdot\hat{\mathbf{r}}=\nabla\phi\cdot\hat{\mathbf{r}}$. Since this gradient can be eliminated by setting $\tilde{\v}=\v-\nabla\phi$, from now on we consider the special case
\[\v=\nabla\times \bs{\psi},\]
where $\v$ is tangent to the boundary. Taking curl of the generalized NS equations \cref{eq_active_fluids} yields the vorticity-vector potential formulation
\begin{subequations}
\label{eq:Ns_psiomega}
\begin{align}
\label{eq:Ns_psiomega_a}
\nabla^2\bs{\psi}&=-\bs{\omega},\\
\label{eq:Ns_psiomega_b}
(\partial_t-\Gamma_0\nabla^2+\Gamma_2\nabla^4-\Gamma_4\nabla^6)\bs{\omega}&=-\nabla\times(\bs{\omega}\times \v),
\end{align}
\end{subequations}
where $\bs{\omega}=\nabla\times \v$ is the vorticity. In the next section, we decouple \cref{eq:Ns_psiomega} by using  a vector decomposition called the poloidal-toroidal decomposition.

\subsection{Reformulation of the NS equations with the PT decomposition}
\label{sec_reform_NS}

A divergence-free vector field $\psi$ in spherical coordinates can be decomposed with the poloidal-toroidal (PT) decomposition into an orthogonal sum of a poloidal and toroidal field as~\cite{backus1986poloidal} $\bs{\psi}=\mathbf{P}+\mathbf{T}$. Moreover, there exist a poloidal and toroidal scalars: $P_{\bs{\psi}}$ and $T_{\bs{\psi}}$, unique up to the addition of an arbitrary function of the radial variable $r\in[0,1]$ such that
\begin{align*}
\mathbf{P} &= \nabla\times\nabla\times(\mathbf{r}P_{\bs{\psi}})=-\mathbf{r}\nabla^2 P_{\bs{\psi}}+\nabla[\partial_r (rP_{\bs{\psi}})],\\
\mathbf{T} &= \nabla\times(\mathbf{r}T_{\bs{\psi}}),
\end{align*}
where $\mathbf{r}:=r\hat{\mathbf{r}}$. Thanks to the above decomposition, we can write down the PT decomposition of the velocity and vorticity fields as
\begin{subequations} \label{eq_decoupled_PT}
\begin{align}
\v &= \nabla\times\psi = \nabla\times\nabla\times(\mathbf{r}T_{\bs{\psi}})+\nabla\times(-\mathbf{r}\nabla^2 P_{\bs{\psi}}), \label{eq:PTv}\\
\bs{\omega} &= \nabla\times \v=\nabla\times\nabla\times(-\mathbf{r}\nabla^2 P_{\bs{\psi}})+\nabla\times(-\mathbf{r}\nabla^2 T_{\bs{\psi}}). \label{eq:PTomega}
\end{align}
\end{subequations}
We have found the action of the poloidal-toroidal decomposition under the transformation $\bs{\psi}\to\nabla\times\nabla\times\bs{\psi}$:
\[P_{\bs{\psi}} \rightarrow -\nabla^2P_{\bs{\psi}},\qquad T_{\bs{\psi}}\rightarrow -\nabla^2 T_{\bs{\psi}}.\]

It follows from \cref{eq:Ns_psiomega_a,eq_decoupled_PT} and the orthogonality of the PT decomposition that the PT scalars for $\bs{\psi}$ and $\bs{\omega}$ obey two Poisson equations
\begin{align*}
\nabla^2 P_{\bs{\psi}}&=-P_{\bs{\omega}}, \\
\nabla^2 T_{\bs{\psi}}&=-T_{\bs{\omega}}.
\end{align*}
Then, the uniqueness of the poloidal and toroidal scalars up to the addition of a function of $r$ and \cref{eq_decoupled_PT} implies that
\begin{subequations}
\label{eq_Poisson_velocity}
\begin{align}
\nabla^2 P_{\v} &= -T_{\bs{\omega}}, \\
T_{\v} &= P_{\bs{\omega}},
\end{align}
\end{subequations}
where made the identifications $P_{\v}=T_{\bs{\psi}}$ and $T_{\v}=-\nabla^2 P_{\bs{\psi}}$. Additionally, $P_{\v}$ satisfies the Dirichlet boundary condition because $\bs \psi$ is normal to the boundary implying that $T_{\bs{\psi}}$ and thus $P_{\v}$ must vanish on the boundary.

We now turn to the PT decomposition of \cref{eq:Ns_psiomega_b} and derive the equations satisfied by the vorticity PT scalars $(P_{\bs{\omega}},T_{\bs{\omega}})$. Again, the orthogonality of the PT decomposition and \cref{eq_decoupled_PT} imply that $(P_{\bs{\omega}},T_{\bs{\omega}})$ satisfy
\begin{subequations}
\label{eq:PTvorticity_eqs}
\begin{align}
(\partial_t-\Gamma_0\nabla^2+\Gamma_2\nabla^4-\Gamma_4\nabla^6)P_{\bs{\omega}} &=-P_{\nabla\times({\bs{\omega}}\times \v)}, \\
(\partial_t-\Gamma_0\nabla^2+\Gamma_2\nabla^4-\Gamma_4\nabla^6)T_{\bs{\omega}} &=-T_{\nabla\times({\bs{\omega}}\times \v)},
\end{align}
\end{subequations}
where the right-hand sides correspond to the PT decomposition of the nonlinear term $\nabla\times(\omega\times v)$. 

\subsection{Boundary conditions}
\label{sec_boundary_conditions}

Since we focus on velocity fields $\v$ that are tangent to the boundary, the most general boundary conditions for $\v$ correspond to an arbitrary two-dimensional vector field $\v^\textnormal{t}$ on the unit sphere. We use the Hodge decomposition to represent $\v^\textnormal{t}$ as a sum of a surface gradient and a surface curl (skew-gradient):
\begin{equation}
\label{eq:Vbc_decomposition}
\v|_{\partial\Omega}=\v^\textnormal{t}=\nabla_1 f+\Lambda_1 g,
\end{equation}
where $f, g$ are two arbitrary functions on the unit sphere. Let $(\hat{\mathbf{r}},\hat{\mathbf{\lambda}},\hat{\mathbf{\theta}})$ be the unit vectors in spherical coordinates in the radial, azimuthal, and polar directions, then the surface gradient and curl are defined as
\[\nabla_1 f = \frac{1}{\sin\theta}\frac{\partial f}{\partial \lambda}\hat{\mathbf{\lambda}}+\frac{\partial f}{\partial\theta}\hat{\mathbf{\theta}},\qquad \Lambda_1 g = -\frac{\partial g}{\partial \theta}\hat{\mathbf{\lambda}}+\frac{1}{\sin\theta}\frac{\partial g}{\partial \lambda}\hat{\mathbf{\theta}}.
\]
The two surface potentials $f$ and $g$ parametrize the compressible and irrotational part of $\v^\textnormal{t}$, respectively, and are user-specified functions, which can be time-dependent.

\subsubsection{Dirichlet conditions for $P_{\bs{\omega}}$}
\label{sec_Dirichlet_P_omega}

We observe that the poloidal scalar $P_{\bs{\omega}}$ obeys the following Dirichlet boundary condition:
\begin{equation}
\label{eq:Pomega_Dirichlet}
P_{\bs{\omega}}|_{\partial\Omega} = g(\lambda,\theta),
\end{equation}
where $g$ is the scalar potential of the irrotational part of the velocity field $\v^\textnormal{t}$, $\lambda\in [-\pi,\pi]$ is the azimuthal variable, and $\theta\in[0,\pi]$ is the polar angle. This follows from \cref{eq_decoupled_PT,eq:Vbc_decomposition}, which imply that $T_{\v}(1,\lambda,\theta)=g(\lambda,\theta)=P_{\bs{\omega}}(1,\lambda,\theta)$. 

\subsubsection{Integral conditions for $T_{\bs{\omega}}$}
\label{sec_Integral_T_omega}

Quartapelle and Valz-Gris~\cite{quartapelle1981projection} showed that specifying the boundary conditions for the velocity field $\v$ is equivalent to imposing integral conditions on the vorticity $\bs{\omega}$. The integral conditions take the form of projections of $\bs{\omega}$ onto the linear space of vector fields $\bs{\eta}$ that satisfy $\nabla\times\nabla\times \bs{\eta}=0$. We first show that only half of these projections correspond to true integral conditions, and the other half is equivalent to specifying the normal component of $\bs{\omega}$.

Consider the linear space of vector fields $\bs{\eta}$ that satisfy $\nabla\times\nabla\times\bs{\eta}=0$. Without loss of generality, we may assume that $\bs{\eta}$ is divergence-free. We apply the PT decomposition to $\bs{\eta}$ and note that $\nabla\times\nabla\times\bs{\eta}=-\nabla^2\bs{\eta}=0$ implies that the corresponding PT scalars are harmonic functions:
\[\nabla^2P_{\bs{\eta}}=0, \qquad \nabla^2 T_{\bs{\eta}}=0,\]
where again we employed the orthogonality of the PT decomposition to obtain two independent equations. Since $P_{\bs{\eta}}$ is harmonic, we immediately conclude that the poloidal fields $\mathbf{P}_{\bs{\eta}}$ are gradients of functions
\[\mathbf{P}_{\bs{\eta}}=\nabla\times\nabla\times (\mathbf{r}P_{\bs{\eta}})=-\mathbf{r}\nabla^2 P_{\bs{\eta}}+\nabla[\partial_r (rP_{\bs{\eta}})]=\nabla[\partial_r (rP_{\bs{\eta}})].\]
Moreover, since $\mathbf{P}_{\bs{\eta}}$ is divergence-free, it is a gradient of a harmonic function. Therefore, projecting $\bs{\omega}$ onto $\mathbf{P}_{\bs{\eta}}$ is equivalent to specifying the gradient part of the Hodge decomposition for $\bs{\omega}$, which in turn is fully determined by the normal component $\bs{\omega}\cdot\hat{\mathbf{r}}$. We conclude that one half of the integral conditions on $\bs{\omega}$ put forward by Quartapelle and Valz-Gris~\cite{quartapelle1981projection} do not constitute true integral conditions in the sense that they are equivalent to boundary conditions. This equivalence is demonstrated by the fact that the poloidal scalar $P_{\bs{\omega}}$ obeys the Dirichlet boundary conditions~\eqref{eq:Pomega_Dirichlet} rather than integral conditions. Physically, this means that the normal component of the vorticity is determined by the surface curl of the velocity field on the boundary. 

We now show that the toroidal scalar $T_{\bs{\omega}}$ obeys integral conditions.
Consider projecting $\bs{\omega}$ onto toroidal fields $\mathbf{T}_{\bs{\eta}}=\nabla\times (\mathbf{r}T_{\bs{\eta}})$, where $T_{\bs{\eta}}$ is a harmonic function. We obtain
\[
\int_{\Omega} \bs{\omega} \cdot \mathbf{T}_{\bs{\eta}}\d x=\int_{\Omega} \mathbf{T}_{\bs{\omega}} \cdot \mathbf{T}_{\bs{\eta}}\d x=\int_{\Omega} \nabla\times(\mathbf{r}T_{\bs{\omega}})\cdot \nabla\times(\mathbf{r}T_{\bs{\eta}})\d x,
\]
where we used the orthogonality of the PT decomposition. Integrating by parts yields
\[\int_{\Omega} \bs{\omega} \cdot \mathbf{T}_{\bs{\eta}}\d x=\int_{\Omega} T_{\bs{\omega}}\mathbf{r}\cdot \nabla\times\nabla\times(\mathbf{r}T_{\bs{\eta}}) \d x = \int_{\Omega} T_{\bs{\omega}}\mathbf{r} \cdot \{-\mathbf{r}\nabla^2 T_\eta+\nabla[\partial_r (rT_\eta)]\}\d x = \int_{\Omega} T_{\bs{\omega}} r \partial_{rr} (rT_{\bs{\eta}})\d x,\]
where we used the fact that $T_{\bs{\eta}}$ is harmonic. We now show that the above true integral conditions are equivalent to specifying the compressible part of the two-dimensional flow fixed by the boundary conditions on the velocity field.

Consider the same integral conditions as above but now integrate by parts to convert $\bs{\omega}$ into $\v$ by using $\bs{\omega}=\nabla\times\v$:
\begin{equation}
 \label{eq:ICs_meaning1}
\int_{\Omega} \bs{\omega} \cdot \mathbf{T}_{\bs{\eta}} \d x 
= \int_{\Omega} \nabla\times\v\cdot \nabla\times(\mathbf{r}T_{\bs{\eta}})\d x
= \int_{\Omega} \nabla\cdot[\v\times \nabla\times(\mathbf{r}_{\bs{\eta}})]\d x + \int_{\Omega} \v\cdot\nabla\times\nabla\times (\mathbf{r}T_{\bs{\eta}}) \d x.
\end{equation}
We first note that the second integral vanishes as
\[
\int_{\Omega} \v\cdot\nabla\times\nabla\times (\mathbf{r}T_{\bs{\eta}}) \d x
= \int_{\Omega} \v\cdot\nabla[\partial_r (rT_{\bs{\eta}})]\d x=0.
\]
This follows from the orthogonality of the Hodge decomposition as the second integral above is the projection of $\v$ onto gradients of functions while $\v$ comes from the vector potential only. Going back to \cref{eq:ICs_meaning1}, we employ the divergence theorem
\begin{align*}
\int_{\Omega} \bs{\omega} \cdot \mathbf{T}_{\bs{\eta}} \d x 
= \int_{\Omega} \nabla\cdot[\v\times \nabla\times(\mathbf{r}T_{\bs{\eta}})] \d x=\int_{\partial \Omega} \hat{\mathbf{r}}\cdot[\v\times \nabla\times(\mathbf{r}T_{\bs{\eta}})] \d s.
\end{align*}
We can use the triple product to rewrite the surface integral as
\[
\hat{\mathbf{r}}\cdot[\v\times \nabla\times(\mathbf{r}T_{\bs{\eta}})]=
-\v\cdot[\hat{\mathbf{r}}\times \nabla\times(\mathbf{r}T_{\bs{\eta}})]=-\v\cdot \nabla_1 T_{\bs{\eta}},
\]
where $\nabla_1$ is the surface gradient. The integral condition reads
\begin{align}
\int_{\Omega} \bs{\omega} \cdot \mathbf{T}_{\bs{\eta}} \d x =
-\int_{\partial\Omega}\v\cdot \nabla_1 T_{\bs{\eta}} \d s =
-\int_{\partial\Omega}\nabla _1f\cdot \nabla_1 T_{\bs{\eta}} \d s=
\int_{\partial\Omega}f\nabla_1^2 T_{\bs{\eta}} \d s,
\end{align}
where in the second equality we used the orthogonality of the Hodge decomposition of velocity fields on the boundary and in the last equality we integrated by parts. We conclude that the true integral conditions, corresponding to the projections of the vorticity onto toroidal fields with harmonic toroidal scalars, fix the compressible part of the velocity field $\v^\textnormal{t}$ on the boundary parameterized by the surface potential $f$ (cf.~\cref{sec_boundary_conditions}).

\section{Numerical method} \label{sec_numer_method}

In this section, we discuss the numerical method used to discretize the poloidal-toroidal formulation of the incompressible NS equations. We then wish to solve the following equations (see \cref{sec_reform_NS}):
\begin{align*}
\nabla^2 P_{\v} &= -T_{\bs{\omega}}, \\
T_{\v} &= P_{\bs{\omega}},
\end{align*}
with homegeneous Dirichlet boundary conditions on $P_{\v}$. The second system of equations reads
\begin{subequations}
\label{eq_Pw_Tw}
\begin{align}
(\partial_t-\Gamma_0\nabla^2+\Gamma_2\nabla^4-\Gamma_4\nabla^6)P_{\bs{\omega}}
&= -P_{\nabla\times(\bs{\omega}\times \v)}, \\
(\partial_t-\Gamma_0\nabla^2+\Gamma_2\nabla^4-\Gamma_4\nabla^6)T_{\bs{\omega}}
&= -T_{\nabla\times(\bs{\omega}\times \v)},
\end{align}
\end{subequations}
with Dirichlet boundary conditions on $P_{\bs{\omega}}$ (see \cref{sec_Dirichlet_P_omega}) and integral conditions on $T_{\bs{\omega}}$ according to \cref{sec_Integral_T_omega}. We discretize \cref{eq_Pw_Tw} in time using the implicit-explicit backward differentiation of order one (IMEX-BDF1) time-stepping scheme. In this paper, we focus on the optimal complexity spatial discretization and employ this time-stepping scheme for simplicity but higher order schemes may be more appropriate for solving the NS equations~\cite{canuto2012spectral,karniadakis1991high,kassam2005fourth,kim1985application}. Let $\Delta t$ be the temporal discretization and $k\geq 0$ the current time-step of the algorithm. At step $k+1$, assuming the scalars $P_\omega^k$ and $T_\omega^k$ have already been computed, we have to solve the following Poisson's equation with zero Dirichlet conditions:
\begin{subequations}
\label{eq_scalar_v}
\begin{align}
\nabla^2 P_{\v}^k &=-T_{\bs{\omega}}^k, \\
T_{\v}^k &= P_{\bs{\omega}}^k.
\end{align}
\end{subequations}
Moreover, if we denote $\mathbf{N}=\nabla\times(\bs{\omega}\times \v)$, then \cref{eq_Pw_Tw} reads
\begin{subequations}
\label{eq_scalar_w}
\begin{align}
(-\Gamma_0\nabla^2+\Gamma_2\nabla^4-\Gamma_4\nabla^6)P_{\bs{\omega}}^{k+1}+\frac{1}{\Delta t}P_{\bs{\omega}}^{k+1}&=\frac{1}{\Delta t}P_{\bs{\omega}}^{k} - P_{\mathbf{N}}^k,\\
(-\Gamma_0\nabla^2+\Gamma_2\nabla^4-\Gamma_4\nabla^6)T_{\bs{\omega}}^{k+1}+\frac{1}{\Delta t}T_{\bs{\omega}}^{k+1}&=\frac{1}{\Delta t}T_{\bs{\omega}}^{k} - T_{\mathbf{N}}^k.
\end{align}
\end{subequations}

We propose a fast algorithm to solve Eqs.~\eqref{eq_scalar_v}-\eqref{eq_scalar_w} and thus achieve an optimal complexity solver to the NS equations per time step. Here, ``optimal'' denotes a linear complexity, up to polylogarithmic factors, in terms of the degrees of freedom needed to resolve the solution spatially. The poloidal and toroidal components of the vorticity and the velocity vector fields are approximated using a global spectral method, introduced in section~\ref{sec_spatial_discretization}. Our choice of spectral basis takes into account the existence of fast --FFT based--  transforms between spectral coefficients and physical values space, which is needed to compute the nonlinear advection term in the NS equations since it is treated explicitly by our time-stepping scheme. This nonlinear vector field also requires a poloidal-toroidal decomposition algorithm and the ability to perform vector calculus operations (see \cref{sec_nonlinear}). Finally, \cref{sec_Helmholtz_SH} consists of the description of a fast Helmholtz's solver with Dirichlet and integral conditions, used to numerically solve \cref{eq_scalar_v,eq_scalar_w} at every time-step. \cref{fig_diagram_NS_solver} summarizes the NS algorithm and shows the resulting optimal complexity per time step.

\tikzstyle{block} = [draw, rectangle, minimum height=3em, minimum width=7em]
\tikzstyle{pinstyle} = [pin edge={to-,thin,black}]
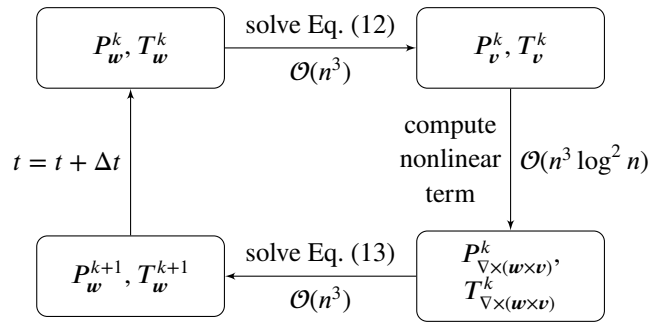
\begin{figure}[htbp]
\centering
\begin{tikzpicture}[auto, node distance=5cm,>=latex']
    \node [block, rounded corners,align=center](Step1) {$P_{\bs{w}}^k$, $T_{\bs{w}}^k$};
    \node [block, rounded corners,align=center, right of=Step1] (Step2) {$P_{\v}^k$, $T_{\v}^k$};
    \node [block, rounded corners,align=center, below of=Step2, node distance=3cm] (Step3) {$P_{\nabla\times (\bs{w}\times \v)}^k$,\\ $T_{\nabla\times (\bs{w}\times \v)}^k$};
    \node [block, rounded corners,align=center, left of=Step3] (Step4) {$P_{\bs{w}}^{k+1}$, $T_{\bs{w}}^{k+1}$};
    \draw [->] (Step1) -- (Step2) node[midway,above] {solve Eq.~\eqref{eq_scalar_v}} node[midway,below] {$\mathcal{O}(n^3)$};
    \draw [->] (Step2) -- (Step3) node[midway,left,align=center] {compute\\nonlinear\\ term} node[midway,right] {$\mathcal{O}(n^3\log^2 n)$};
    \draw [->] (Step3) -- (Step4) node[midway,above] {solve Eq.~\eqref{eq_scalar_w}} node[midway,below] {$\mathcal{O}(n^3)$};
    \draw [->] (Step4) -- (Step1) node[midway,left] {$t = t+\Delta t$};
\end{tikzpicture}
\caption{Tasks performed during one time-step of the NS solver with associated computational complexity. The poloidal and toroidal scalars $P_{\bs{w}}$, $T_{\bs{w}}$ are discretized using $\mathcal{O}(n^3)$ physical values ($n$ for each radial, azimuthal, and polar direction).}
\label{fig_diagram_NS_solver}
\end{figure}

\subsection{Spatial discretization of functions and vectors on the ball}
\label{sec_spatial_discretization}

We now discuss the different polynomial basis employed to approximate functions and vector-valued functions defined on the ball, perform vector calculus operations, compute poloidal-toroidal decompositions and solve Helmholtz's equation.

\subsubsection{Chebyshev--Fourier--Fourier basis}
\label{sec_CFF_basis}

If $f$ is a smooth function on the ball, then it can be represented in a Chebyshev--Fourier--Fourier (CFF) series:
\begin{equation}
\label{eq_CFF_series_inf}
f(r,\lambda,\theta)\approx \sum_{k=0}^{+\infty}\sum_{l=-\infty}^{+\infty}\sum_{m=-\infty}^{+\infty}f_{klm}T_k(r)e^{\mathbf{i}l\lambda}e^{\mathbf{i}m\theta},
\end{equation}
where $(r,\lambda,\theta)\in[0,1]\times[-\pi,\pi]\times[0,\pi]$ and $T_k$ is the degree $k$ Chebyshev polynomial of the first kind. The function $f$ is discretized in $\mathcal{O}(n^3)$ coefficients by truncating the series~\eqref{eq_CFF_series_inf},
\begin{equation}
\label{eq_CFF_series}
f(r,\lambda,\theta)\approx \sum_{k=0}^{n/2}\sum_{l=-n/2}^{n/2}\sum_{m=-n/2}^{n/2}f_{klm}T_k(r)e^{\mathbf{i}l\lambda}e^{\mathbf{i}m\theta},
\end{equation}
where $n$ is an even integer.

A few salient features make this discretization choice attractive for working with functions on the ball. First, any smooth function can be represented as in~\eqref{eq_CFF_series} using the double sphere method~\cite{merilees1973pseudospectral} to remove artificial boundary conditions at the poles and at the origin. This method, originally proposed by Merilees, has been extended to various tensor-product domains such as the sphere, the disk, and the ball~\cite{boulle2020computing,townsend2016computing,wilber2017computing}. This method extends a function $f$ defined on the ball to a function $\tilde{f}$ defined on $[-1,1]\times[-\pi,\pi]\times[-\pi,\pi]$. Moreover, the Chebyshev--Fourier--Fourier coefficients of $\tilde{f}$, $(f_{klm})$, can be computed in $\mathcal{O}(n^3\log n)$ operations by sampling $f$ at the three dimensional Chebyshev--Fourier--Fourier grid
\begin{equation}
\label{eq_CFF_grid}
\left(\cos\left(\frac{2k\pi}{n}\right),\frac{2l\pi}{n},\frac{2m\pi}{n}\right),\qquad 0\leq k\leq \frac{n}{2},\quad -\frac{n}{2}\leq l,m\leq \frac{n}{2},
\end{equation}
and using transforms based on the FFT~\cite{cooley1965algorithm}. Finally, operations such as differentiation, vector calculus, and computing the poloidal-toroidal decomposition of vector-valued functions can be done efficiently and are already in the package Ballfun~\cite{boulle2020computing}, which is a component of the Chebfun software~\cite{Driscoll2014}.

\subsubsection{Chebyshev-Spherical harmonic basis}
\label{sec_CSH_basis}

One can also approximate a smooth function $f$ on the ball by the following expansion, called a Chebyshev-Spherical harmonic (CSH) series,
\begin{equation}
\label{eq_CSH_inf}
f(r,\lambda,\theta) \approx \sum_{k=0}^{+\infty}\sum_{l=0}^{+\infty}\sum_{m=-l}^l f_{klm}T_k(r)Y^m_l(\theta,\lambda).
\end{equation}
The spherical harmonics $Y^m_l$ are defined by $Y^m_l(\theta,\lambda) = \tilde{P}^m_l(\cos\theta)e^{im\lambda}$, where $P^m_l$ is the normalized associated Legendre polynomials of degree $l$ and order $m$ such that
\[\int_{\partial\Omega}Y^m_l Y^{m'*}_{l'} \d s = \delta_{ll'}\delta_{mm'},\]
where $\delta_{ll'}=1$ if $l=l'$ and zero otherwise. The truncation of~\eqref{eq_CSH_inf} to Chebyshev polynomials up to degree $n/2$ and wave frequencies up to order $n/2$ yields
\begin{equation}
\label{eq_CSH}
f(r,\lambda,\theta) \approx \sum_{k=0}^{n/2}\sum_{l=0}^{n/2}\sum_{m=-l}^l f_{klm}T_k(r)Y^m_l(\theta,\lambda).
\end{equation}

Slevinsky introduced a fast and practical algorithm, called the spherical harmonic transform (SHT), for transformations between spherical harmonic expansions and bivariate Fourier series~\cite{slevinsky2017fast}. We use this algorithm to convert functions expressed in the CSH basis into the CFF basis in $\mathcal{O}(n^3\log^2n)$ operations. The inverse discrete cosine transform (DCT)~\cite{gentleman1972implementing,gentleman1972implementing2,mason2002chebyshev} and the fast Fourier transform~\cite{cooley1965algorithm} are then used to convert the CFF coefficients into values. The following diagram summarizes the operations used to evaluate a function expressed as a CSH series at the CFF points:
\[\text{CSH coefficients} \xlongleftrightarrow[\mathcal{O}(n^3\log^2 n)]{\text{DCT+SHT}} \text{CFF coefficients} \xlongleftrightarrow[\mathcal{O}(n^3\log n)]{\text{DCT+FFT}} \text{Values at CFF points}.\]

\subsection{Vector calculus and PT decomposition}
\label{sec_nonlinear}

In this section, we assume the existence of two functions $P_{\v}$ and $T_{\v}$, representing the poloidal and toroidal scalars (see \cref{sec_reform_NS}) of a divergence-free vector field $\v$ in the CSH basis. These functions are approximated by a series of $\mathcal{O}(n^3)$ Chebyshev and spherical harmonics coefficients and the number of degrees of freedom needed to represent the three components of the velocity field $\v$ is equal to $N=3(n/2+1)(n+1)^2=\mathcal{O}(n^3)$. We then explain the computation the poloidal-toroidal decomposition of $\nabla\times (\bs{\omega}\times \v)$, where $\bs{\omega}$ is a vector-valued function defined by $\bs{\omega}=\nabla\times \v$.

Given the functions $P_{\v}$ and $T_{\v}$ in the CSH basis, we first convert them to the CFF basis with the spherical harmonic transform (see \cref{sec_CSH_basis}) in $\mathcal{O}(n^3\log^2n)$ operations. According to \cref{sec_reform_NS}, the vector-valued function $\v$ can be recovered from its poloidal and toroidal scalars using vector calculus:
\[\v=\nabla\times\nabla\times(\mathbf{r}P_{\v}) + \nabla\times(\mathbf{r}T_{\v}).\]
This last operation is performed in optimal complexity, $\mathcal{O}(n^3)$, in the CFF basis~\cite{boulle2020computing}. The vector field $\bs{\omega} = \nabla\times \v$ is computed as well in $\mathcal{O}(n^3)$ operations using the expression of $\v$ in the CFF basis. This allows us to perform only two forward SHTs (on the poloidal and toroidal scalars of $\v$) to compute $\v$ and $\bs{\omega}$ in the CFF basis. We then sample the Cartesian components of $\v$ and $\bs{\omega}$ with six fast transforms based on FFT in $\mathcal{O}(n^3\log n)$ operations (see \cref{sec_CFF_basis}). The tensors representing the evaluation of $\v$ and $\bs{\omega}$ at the CFF points are multiplied element-wise and transformed back into CFF coefficients to get an approximation of $\bs{\omega}\times \v$ as a CFF series. We then take the curl of $\bs{\omega}\times \v$ by doing vector calculus operations in the CFF basis. The next step is to compute the poloidal-toroidal decomposition of $\nabla\times \bs{\omega}\times \v$. A fast algorithm in $\mathcal{O}(n^3)$ operations is described in~\cite{boulle2020computing} and returns the CFF coefficients of the poloidal and toroidal scalars $P_{\nabla\times(\bs{\omega}\times \v)}$ and $T_{\nabla\times(\bs{\omega}\times \v)}$. Finally, these functions, expressed as CFF series, are transformed in the CSH basis via an inverse SHT. \cref{fig_diagram_nonlinear} summarizes the algorithms used to compute the nonlinear term. Note that we expand the CFF series representing the Cartesian components of $w$ and $v$ using the 3/2 rule to prevent aliasing effects during the computation of $\omega\times v$~\cite{kirby2003aliasing}. Later, we truncate the nonlinear term to the original spatial discretization.

\tikzstyle{block} = [draw, rectangle, minimum height=3em, minimum width=7em]
\tikzstyle{pinstyle} = [pin edge={to-,thin,black}]
\begin{figure}[htbp]
\centering
\begin{tikzpicture}[auto, node distance=5cm,>=latex']
    \node [block, rounded corners,align=center](Step1) {$P_{\v}$, $T_{\v}$\\ in CSH basis};
    \node [block, rounded corners,align=center, right of=Step1] (Step2) {$P_{\v}$, $T_{\v}$\\ in CFF basis};
    \node [block, rounded corners,align=center, right of=Step2] (Step3) {$\v$, $\bs{\omega}$};
    \node [block, rounded corners,align=center, below of=Step3, node distance=3cm] (Step4) {$\tilde{\v}$, $\tilde{\bs{\omega}}$\\ in values space};
    \node [block, rounded corners,align=center, left of=Step4] (Step5) {$\widetilde{\bs{\omega}\times \v}$\\ in values space};
    \node [block, rounded corners,align=center, left of=Step5] (Step6) {$\bs{\omega}\times \v$};
    \node [block, rounded corners,align=center, below of=Step6, node distance=3cm] (Step7) {$\nabla\times (\bs{\omega}\times \v)$};
    \node [block, rounded corners,align=center, right of=Step7] (Step8) {$P_{\mathbf{N}}$, $T_{\mathbf{N}}$\\ in CFF basis};
    \node [block, rounded corners,align=center, right of=Step8] (Step9) {$P_{\mathbf{N}}$, $T_{\mathbf{N}}$\\ in CSH basis};

    \draw [->] (Step1) -- (Step2) node[midway,above] {DCT, SHT} node[midway,below] {$\mathcal{O}(n^3\log^2 n)$};
    \draw [->] (Step2) -- (Step3) node[midway,above] {vector calculus} node[midway,below] {$\mathcal{O}(n^3)$};
    \draw [->] (Step3) -- (Step4) node[midway,left,align=center] {$\mathcal{O}(n^3\log n)$} node[midway,right,align=center] {DCT,\\ FFT};
    \draw [->] (Step4) -- (Step5) node[midway,above] {multiplication} node[midway,below] {$\mathcal{O}(n^3)$};
    \draw [->] (Step5) -- (Step6) node[midway,above] {DCT, FFT} node[midway,below] {$\mathcal{O}(n^3\log n)$};
    \draw [->] (Step6) -- (Step7) node[midway,left,align=center] {$\mathcal{O}(n^3)$} node[midway,right,align=center] {vector\\ calculus};
    \draw [->] (Step7) -- (Step8) node[midway,above, align=center] {PT\\ decomposition} node[midway,below] {$\mathcal{O}(n^3)$};
    \draw [->] (Step8) -- (Step9) node[midway,above] {DCT, SHT} node[midway,below] {$\mathcal{O}(n^3\log^2 n)$};
\end{tikzpicture}
\caption{Procedure used to compute the poloidal-toroidal decomposition of the vector field $\mathbf{N}:=\nabla\times (\bs{\omega}\times \v)$ in the CSH basis. $\tilde{\v}$, $\tilde{\bs{\omega}}$, and $\widetilde{\bs{\omega}\times \v}$ respectively denote the CFF values of the vector-valued functions $\v$, $\bs{\omega}$, and $\bs{\omega}\times \v$.}
\label{fig_diagram_nonlinear}
\end{figure}
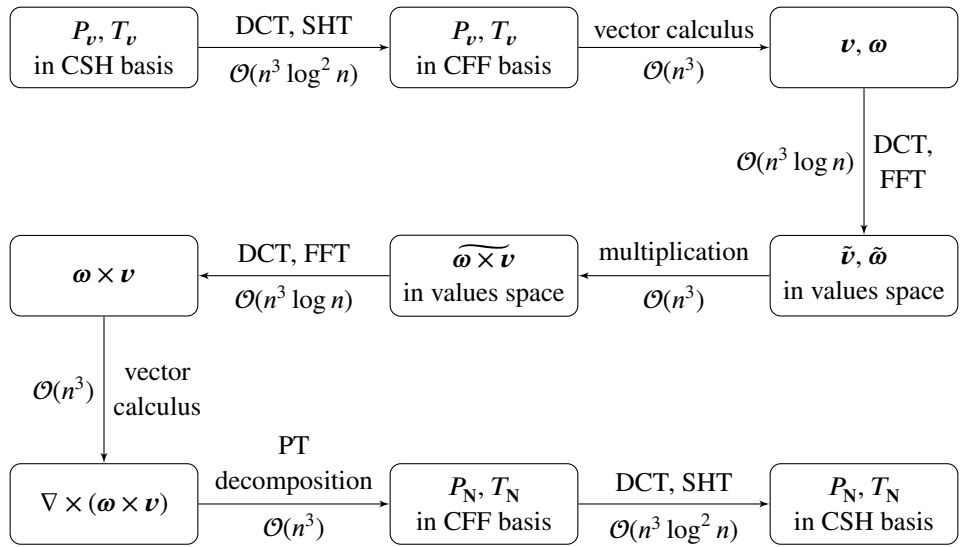

\subsection{Helmholtz's solver}
\label{sec_Helmholtz_SH}

This section describes an optimal complexity algorithm for solving Helmholtz's equation on the ball with Dirichlet or integral boundary conditions. Helmholtz's equation with Dirichlet boundary conditions, discretized using the ultraspherical spectral method~\cite{olver2013fast} in the Chebyshev--Fourier--Fourier basis, yields to a Sylvester matrix equation, which can be solved in optimal complexity~\cite{fortunato2017fast}. However, it is not clear whether solving Helmholtz's equation with integral conditions in optimal complexity is possible in this basis. We then choose to represent functions and solve Helmholtz's equation in the Chebyshev-Spherical harmonic basis since spherical harmonics decouple the angular part of the Laplacian in spherical coordinates. This choice, associated with the ultraspherical spectral method, introduced in~\cite{olver2013fast}, leads to a fast and optimal complexity Helmholtz solver. Note that the algorithm described in this section is also applicable to ``higher-order'' Helmholtz equations such as  
Eqs.~\eqref{eq_scalar_v}-\eqref{eq_scalar_w} but the description of the numerical method becomes more complicated.

\subsubsection{Spatial discretization of Helmholtz's equation}
\label{sec_disc_helmholtz}

Consider the following Helmholtz's equation on the ball
\begin{equation}
\label{eq_helmholtz_spherical}
\nabla^2u + K^2u = f,
\end{equation}
where $\nabla^2$ stands for the laplacian and $K$ is the wave number. We first begin by writing \cref{eq_helmholtz_spherical} in spherical coordinate using the change of variables $(x,y,z)=(r\cos\lambda\sin\theta,r\sin\lambda\sin\theta,r\cos\theta)$:
\begin{equation}
\label{eq_helmholtz_spherical_coordinates}
\frac{1}{r^2}\frac{\partial}{\partial r}\left(r^2\frac{\partial u}{\partial r}\right) + \frac{1}{r^2\sin\theta}\frac{\partial}{\partial\theta}\left(\sin\theta\frac{\partial u}{\partial\theta}\right)+\frac{1}{r^2\sin^2\theta}\frac{\partial^2u}{\partial\lambda^2} + K^2 u = f,
\end{equation}
where $r\in[-1,1]$, $\lambda\in[-\pi,\pi]$, and $\theta\in[0,\pi]$. Note that the radial variable $r$ has been ``doubled-up'' and extended from $[0,1]$ to $[-1,1]$ to remove the artificial boundary condition at $r=0$ (see \cref{sec_CFF_basis}). \cref{eq_helmholtz_spherical_coordinates} is then multiplied by $r^2$ to remove the singularity at the origin:
\begin{equation}
\label{eq_helmholtz_spherical_2}
\frac{\partial}{\partial r}\left(r^2\frac{\partial u}{\partial r}\right) + \frac{1}{\sin\theta}\frac{\partial}{\partial\theta}\left(\sin\theta\frac{\partial u}{\partial\theta}\right)+\frac{1}{\sin^2\theta}\frac{\partial^2u}{\partial\lambda^2} + K^2 r^2u = r^2 f.
\end{equation}
Furthermore, we seek for a solution $u$ expressed as the following CSH series
\begin{equation}
\label{eq_u_CSH_finite}
u(r,\theta,\lambda) \approx \sum_{k=0}^{n/2}\sum_{l=0}^{n/2}\sum_{m=-l}^l u_{klm}T_k(r)Y^m_l(\theta,\lambda).
\end{equation}
Moreover, for $0\leq l\leq n/2$ and $-l\leq m\leq l$, the spherical harmonic function of degree $l$ and order $m$, $Y^m_l$, is an eigenvector of the surface laplacian $\nabla^2_1$ with corresponding eigenvalue $-l(l+1)$ since
\begin{equation}
\label{eq_SH_eigenvector}
\nabla^2_1Y^m_l:=\frac{1}{\sin\theta}\frac{\partial}{\partial\theta}\left(\sin\theta\frac{\partial Y^m_l}{\partial\theta}\right)+\frac{1}{\sin^2\theta}\frac{\partial^2Y^m_l}{\partial\lambda^2}=-l(l+1)Y^m_l.
\end{equation}
Then, \cref{eq_helmholtz_spherical_2} is equivalent to
\begin{equation}
\label{eq_helmholtz_coupled}
\sum_{k=0}^{n/2}\sum_{l=0}^{n/2}\sum_{m=-l}^l u_{klm}\left[\frac{\partial}{\partial r}\left(r^2\frac{\partial T_k(r)}{\partial r}\right)-l(l+1)T_k(r)+K^2r^2T_k(r)\right]Y^m_l=r^2f.
\end{equation}
The expression of the right hand side of \cref{eq_helmholtz_coupled} in the CSH basis reads
\[r^2f(r,\lambda,\theta)=\sum_{k=0}^{n/2}\sum_{l=0}^{n/2}\sum_{m=-l}^lf_{klm}T_k(r)Y^m_l(\theta,\lambda).\]
We now fix $0\leq l\leq n/2$, $-l\leq m\leq l$, and define the functions $u^l_m$, $f^l_m$ on $r\in[-1,1]$ by
\[u^m_l(r)=\sum_{k=0}^{n/2}u_{klm}T_k(r),\qquad f^m_l(r)=\sum_{k=0}^{n/2}f_{klm}T_k(r).\]
The orthogonality of the spherical harmonics basis allows us to decouple \cref{eq_helmholtz_coupled} into the following second order ordinary differential equation:
\begin{equation}
\label{eq_Helmholtz_spherical_decoupled}
\frac{\partial}{\partial r}\left(r^2\frac{\partial u^m_l(r)}{\partial r}\right)+(K^2r^2-l(l+1))u^m_l(r)=f^m_l(r),\quad r\in[-1,1].
\end{equation}
This equation is solved using the ultraspherical spectral method~\cite{olver2013fast} as it typically leads to almost banded linear systems and optimal complexity numerical solvers. Let $X^m_l$ (resp.~$F^m_l$) be the vector of Chebyshev coefficients $(u_{klm})_{0\leq k\leq n/2}$ (resp.~$(f_{klm})_{0\leq k\leq n/2}$) and $A_l$ be the sparse and banded matrix representing the operator $u\mapsto\frac{\partial}{\partial r}\left(r^2\frac{\partial u}{\partial r}\right)+(K^2r^2-l(l+1))u$ from the Chebyshev basis $T$ to the ultraspherical $C^{(2)}$ basis. This transforms \cref{eq_Helmholtz_spherical_decoupled} into a linear system with an unknown column of Chebyshev coefficients $X^m_l$:
\begin{equation} \label{eq_Helmholtz_spherical_linear}
A_lX^m_l=F^m_l.
\end{equation}

\begin{figure}[ht]
\centering
\begin{overpic}[width=.4\textwidth,trim={40 0 40 0},clip]{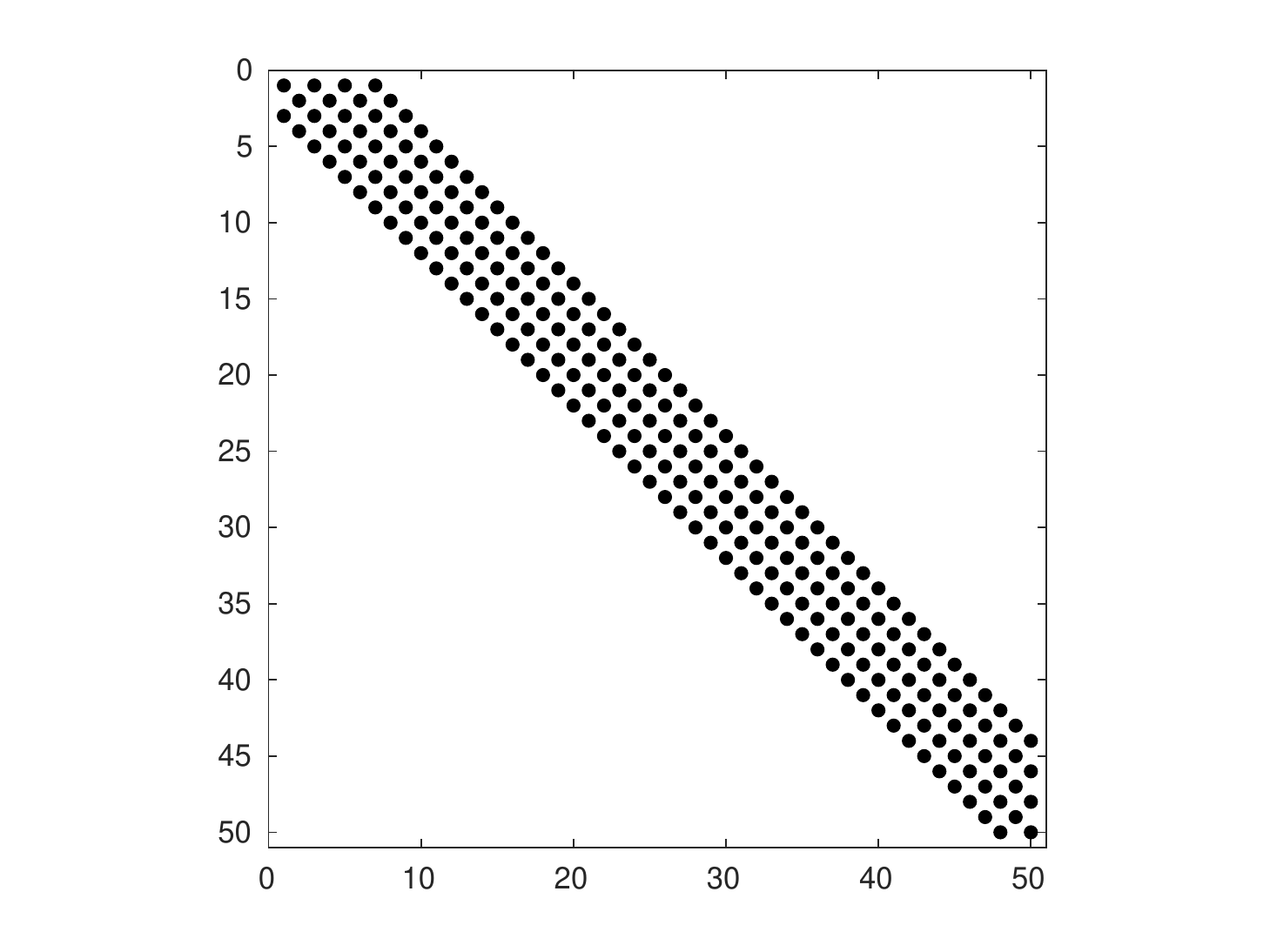}
\end{overpic}
\hspace{1cm}
\centering
\begin{overpic}[width=.4\textwidth,trim={40 0 40 0},clip]{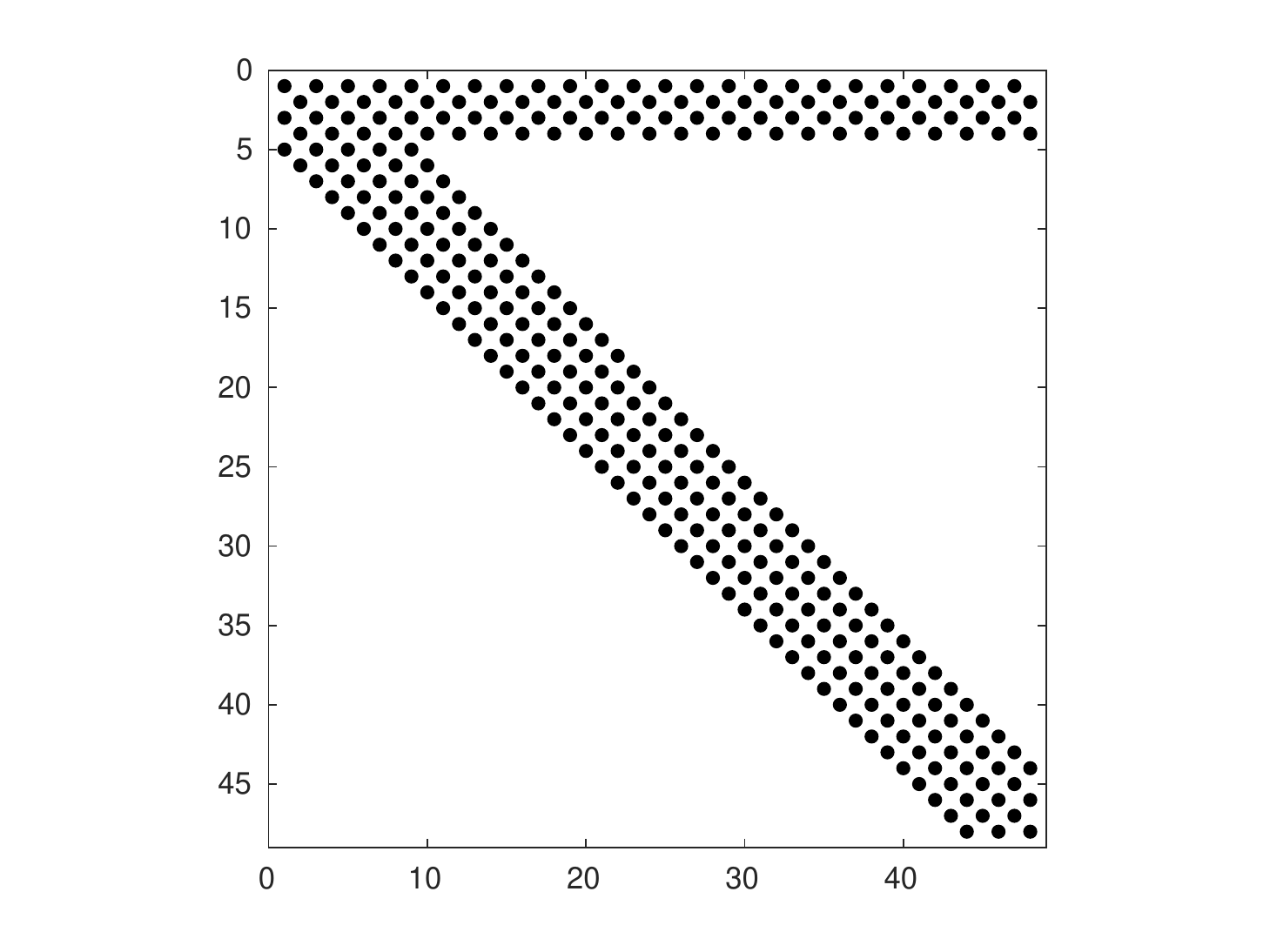}
\end{overpic} 
\caption{Left: Sparsity structure of the matrix $A_l$ for $0\leq l\leq n$. Right: Sparsity structure of the almost banded linear system after adding the boundary conditions.}
\label{fig_sparse_helmholtz}
\end{figure}

\cref{fig_sparse_helmholtz} (left) shows the sparsity pattern of the matrix $A_l$, which is a sparse and banded matrix with bandwidth two. However, $A_l$ is not a full rank matrix because the linear system \eqref{eq_Helmholtz_spherical_linear} does not have a unique solution as \cref{eq_Helmholtz_spherical_decoupled} is a second-order ODE. We then have to modify the matrix $A_l$ to remove degrees of freedom using Dirichlet or integral conditions (see \cref{sec_Dirichlet,sec_numerical_integral}). 

\subsubsection{Dirichlet conditions}
\label{sec_Dirichlet}

We now consider Helmholtz's equation~\eqref{eq_helmholtz_spherical} with Dirichlet boundary conditions
\[u|_{\partial \Omega} = g,\]
where $g$ is a smooth function defined on the unit sphere. $g$ can be represented as a spherical harmonics series
\[g=\sum_{l=0}^{n/2}\sum_{m=-l}^l g_{lm}Y^m_l.\]
The boundary conditions decouple as well as the Helmholtz's equation in the CSH basis (see \cref{sec_disc_helmholtz}) and we obtain the following well-posed problem,
\begin{subequations}
\label{eq_helmholtz_ODE}
\begin{align}
\frac{\partial}{\partial r}\left(r^2\frac{\partial u^m_l}{\partial r}\right)+(K^2r^2-l(l+1))u^m_l=f^m_l, \label{subeq_helmholtz}\\
u_l^m|_{r=1} = g_{lm},\qquad
u_l^m|_{r=-1} = (-1)^lg_{lm}.\label{subeq_helmholtz_2}
\end{align}
\end{subequations}
\cref{subeq_helmholtz} is expressed as a ill-posed linear system (see \cref{eq_Helmholtz_spherical_linear}) in the unknown $X^m_l$, which represents the Chebyshev coefficients of the one-dimensional function $u^m_l$. Likewise, we express \cref{subeq_helmholtz_2} as a linear condition on $X^m_l$,
\begin{equation}
\label{eq_Dirichlet_BC_Helmholtz}
\begin{pmatrix}
1 &  1 & 1 &  1 &\cdots\\
1 & -1 & 1 & -1 &\dots
\end{pmatrix}
X^m_l = 
\begin{pmatrix}
g_{lm}\\
(-1)^l g_{lm}
\end{pmatrix}.
\end{equation}
Incorporating the Dirichlet boundary conditions~\eqref{eq_Dirichlet_BC_Helmholtz} into \cref{eq_Helmholtz_spherical_linear} changes it to a sparse and almost banded linear system, which can be solved in $\mathcal{O}(n)$ operations using an algorithm based on a QR factorization~\cite{olver2013fast}. The resulting linear system's sparsity pattern is shown in \cref{fig_sparse_helmholtz} (right). Finally, solving Helmholtz's equation with Dirichlet boundary conditions requires $\mathcal{O}(n^3)$ operations.

\begin{figure}[ht]
\centering
\begin{overpic}[width=0.4\textwidth,trim={130 60 50 70},clip]{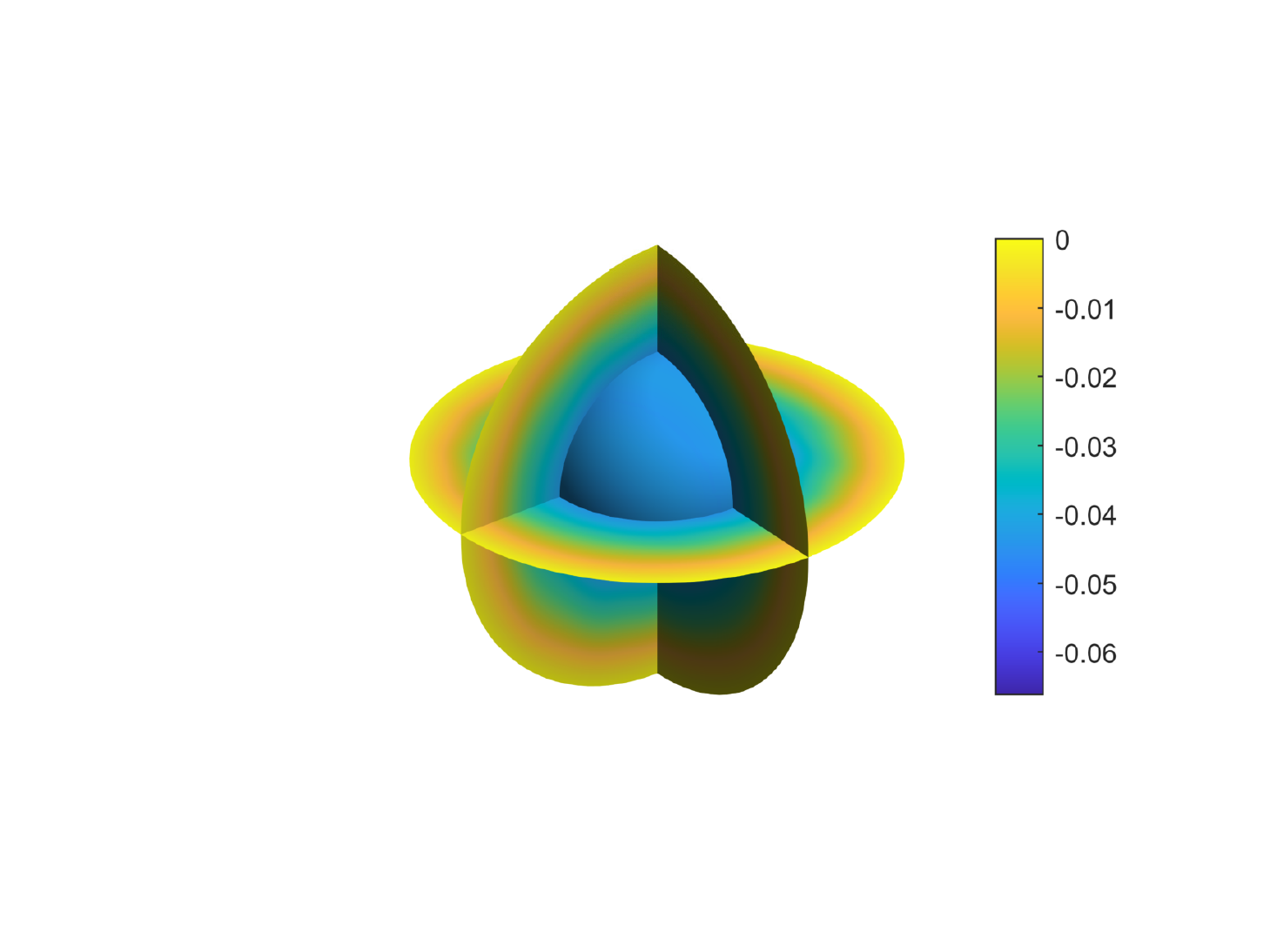}
\end{overpic}
\hspace{2cm}
\begin{overpic}[width=0.36\textwidth]{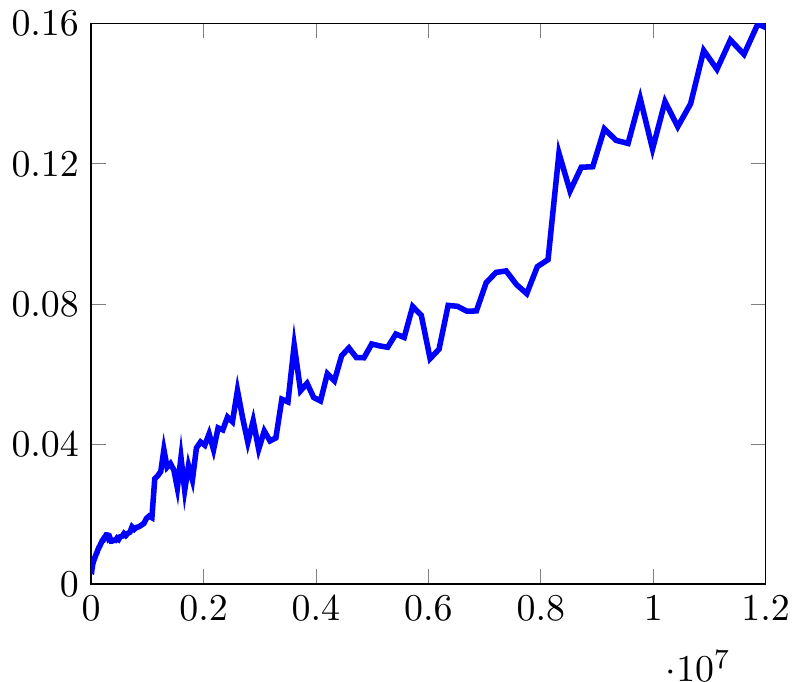}
\put(-7,27.5){\rotatebox{90}{Execution time (s)}}
\put(31,0){Degrees of Freedom}
\end{overpic}
\caption{Left: Solution to the Helmholtz's equation $\nabla^2 u+u = (1-(x^2+y^2+z^2))\cos(10x+5y)^2$ with homegeneous Dirichlet boundary conditions. Right: Execution time for the Helmholtz solver in the CSH basis with Dirichlet boundary conditions.}
\label{fig_helmholtz_example}
\end{figure}

\cref{fig_helmholtz_example} shows a computed solution to the Helmholtz's equation $\nabla^2u+u = (1-(x^2+y^2+z^2))\cos(10x+5y)^2$ with zero Dirichlet boundary conditions and confirms the optimal complexity of the solver in the CSH basis.

\subsubsection{Integral conditions}
\label{sec_numerical_integral}

In this section, we consider Helmholtz's equation of the form
\begin{equation}
\label{eq_Helmholtz_integral}
\nabla^2u+K^2u=f
\end{equation}
with the integral conditions defined in \cref{sec_Integral_T_omega}. After explaining the discretization of the integral conditions in the Chebyshev-Spherical harmonic basis, we describe an optimal complexity algorithm to solve \cref{eq_Helmholtz_integral}. 

The integral conditions read (see \cref{sec_Integral_T_omega}):
\begin{equation}
\label{eq_int_conditions}
\int_{\Omega} ru\partial_{rr}(rT_{\bs{\eta}})^*\d x=-\int_{\partial\Omega}\nabla_1 p\cdot \nabla_1 T_{\bs{\eta}}^* \d s,
\end{equation}
where $p$ is a surface potential defined in \cref{sec_boundary_conditions} and $T_{\bs{\eta}}$ is any harmonic function on the ball. Then, $T_{\bs{\eta}}$ satisfies
\begin{equation} \label{eq_lambda_eta}
\nabla^2 T_{\bs{\eta}}=0,\qquad T_{\bs{\eta}}|_{\partial \Omega}=\alpha,
\end{equation}
where $\alpha$ is an arbitrary smooth function on the sphere. We write $\alpha$ as a series of spherical harmonics
\[\alpha(\lambda,\theta)=\sum_{l=0}^{+\infty}\sum_{m=-l}^l \alpha_{ml} Y^m_l(\lambda,\theta),\]
which transforms the integral conditions into
\begin{equation}
\label{eq_int_conditions_2}
\int_{\Omega} ru\partial_{rr}(rT_{\bs{\eta}}^{m,l})^*\d x=-\int_{\partial\Omega}\nabla_1 p\cdot \nabla_1 T_{\bs{\eta}}^{m,l*} \d s,\qquad l\in\mathbb{N},\quad-l\leq m\leq l,
\end{equation}
where $T_\eta^{m,l}$ denotes the solution to Laplace's equation with Dirichlet boundary conditions $Y^m_l$, i.e.,
\begin{equation} \label{eq_lambda_eta_2}
\nabla^2 T_{\bs{\eta}}^{m,l}=0,\qquad T_{\bs{\eta}}^{m,l}|_{\partial \Omega}=Y^m_l.
\end{equation}
Laplace's equation in spherical coordinates reads
\begin{equation}
\label{eq_laplace_spherical}
\frac{1}{r^2}\frac{\partial}{\partial r}\left(r^2\frac{\partial T_{\bs{\eta}}^{m,l}}{\partial r}\right)+\frac{1}{r^2\sin\theta}\frac{\partial}{\partial\theta}\left(\sin\theta\frac{\partial T_{\bs{\eta}}^{m,l}}{\partial\theta}\right)+\frac{1}{r^2\sin^2\theta}\frac{\partial^2 T_{\bs{\eta}}^{m,l}}{\partial\lambda^2}=0.
\end{equation}
Since spherical harmonics are eigenvectors of the spherical Laplacian $\nabla^2_1$ (see \cref{eq_SH_eigenvector}), they decouple \cref{eq_lambda_eta_2}. Hence, $T_{\bs{\eta}}^{m,l}$ can be written as $T_{\bs{\eta}}^{m,l} = R^m_l Y^m_l$, where $R^m_l$ is a function of the radial variable $r$ to be determined by solving \cref{eq_laplace_spherical}. Then, $R^m_l$ is a solution to the following ordinary differential equation,
\[\frac{\partial}{\partial r}\left(r^2\frac{\partial R^m_l}{\partial r}\right)-l(l+1)R^m_l=0,\]
whose solutions can be expressed as $R^m_l(r)=Ar^{l}+Br^{-(l+1)}$ for $A,B$ are real constants. According to \cref{eq_lambda_eta_2}, $R^m_l$ is a smooth function on the ball and $R^m_l(1)=1$, which implies $A=1$ and $B=0$. Hence, the solution to \cref{eq_lambda_eta_2} is $T_{\bs{\eta}}^{m,l} = r^l Y^m_l$.

Finally, after integrating by parts the right hand side of \cref{eq_int_conditions}, the integral conditions  become
\begin{equation}
\label{eq_int_conditions_3}
\int_{\Omega} ru\partial_{rr}(rT_{\bs{\eta}}^{m,l})^*\d x=\int_{\partial\Omega}p\nabla_1^2 T_{\bs{\eta}}^{m,l*} \d s,\qquad l\in\mathbb{N},\quad -l\leq m\leq l.
\end{equation}
\cref{eq_int_conditions_3} is discretized using the expression of $u$ in the CSH basis and the spherical harmonics expansion of $p$,
\begin{equation}
\label{eq_uf_CSH}
u = \sum_{k=0}^{+\infty}\sum_{l=0}^{+\infty}\sum_{m=-l}^l u_{klm}T_kY^m_l,\quad
f = \sum_{l=0}^{+\infty}\sum_{m=-l}^l p_{lm}Y^m_l.
\end{equation}
Moreover, thanks to the orthonormality of the spherical harmonics basis, the left hand side of \cref{eq_int_conditions_3} is equivalent to
\begin{equation}
\label{eq_integral_LHS}
\int_{\Omega} ru\partial_{rr}(rT_{\bs{\eta}}^{m,l})^*\d x =
\int_0^1 r\sum_{k=0}^{+\infty}u_{klm}T_k(r)\partial_{rr}(r^{l+1})r^2 \d r =l(l+1)\sum_{k=0}^{+\infty}\left(\int_0^1 r^{l+2}T_k(r)\d r\right)u_{klm},
\end{equation}
for every integer $l$ and $-l\leq m\leq l$. Likewise, the right hand side of \cref{eq_int_conditions_3} leads us to the following equality,
\begin{equation}
\label{eq_integral_RHS}
\int_{\partial\Omega}p\nabla_1^2 T_{\bs{\eta}}^{m,l*} \d s = -l(l+1)p_{lm},
\end{equation}
since the spherical harmonic $Y^m_l$ is an eigenvector of the surface Laplacian $\nabla_1^2$ with eigenvalue $-l(l+1)$. By combining \cref{eq_integral_LHS,eq_integral_RHS}, we obtain the discretized formulation of the integral conditions:
\[\sum_{k=0}^{+\infty}\left(\int_0^1 r^{l+2}T_k(r)\d r\right)u_{klm} = -p_{lm},\qquad l>0,\quad -l\leq m\leq l.\]
However, $u$ is written as a finite CSH series (see \cref{eq_u_CSH_finite}), which leads to the following finite formulation of the integral conditions,
\begin{equation}
\label{eq_finite_integral}
\sum_{k=0}^{n/2}\left(\int_0^1 r^{l+2}T_k(r)\d r\right)u_{klm} = -p_{lm},\qquad 1\leq l\leq n,\quad -l\leq m\leq l.
\end{equation}
\cref{eq_finite_integral} is not defined for $l=0$ because the solution to Helmholtz's equation with integral conditions is unique up to the addition of an arbitrary function of $r$. Hence, we impose $u_{k00}=0$ for $0\leq k\leq n/2$ to ensure the uniqueness of the solution to \cref{eq_Helmholtz_integral}.

The integrals in \cref{eq_finite_integral} are computed via Clenshaw--Curtis quadrature by evaluating the Chebyshev polynomials $T_0,...,T_n$ at Chebyshev nodes using Clenshaw's algorithm~\cite{trefethen2013approximation}. This allows us to reformulate \cref{eq_finite_integral} as a linear condition on $X^m_l=(u_{klm})_k$,
\begin{equation}
\label{eq_integral_linear}
B_l X^m_l = -p_{lm}.
\end{equation}
According to \cref{sec_Dirichlet}, this condition can be enforced in the linear system \cref{eq_Helmholtz_spherical_linear}:
\[A_l X^m_l=F^m_l,\]
by replacing the last row of $A_l$ by the integral row condition $B_l$. Then, a permutation of the last row with the first row makes it close to upper triangular (see \cref{fig_sparse_helmholtz}). Finally, Helmholtz's equation with integral conditions is solved in $\mathcal{O}(n^3)$ operations in the CSH basis as inverting this linear system demands $\mathcal{O}(n)$ operations for $1\leq l\leq n/2$ and $-l\leq m\leq l$.

\section{Numerical examples} \label{sec_numer_examples}

\subsection{Convergence test and execution time}

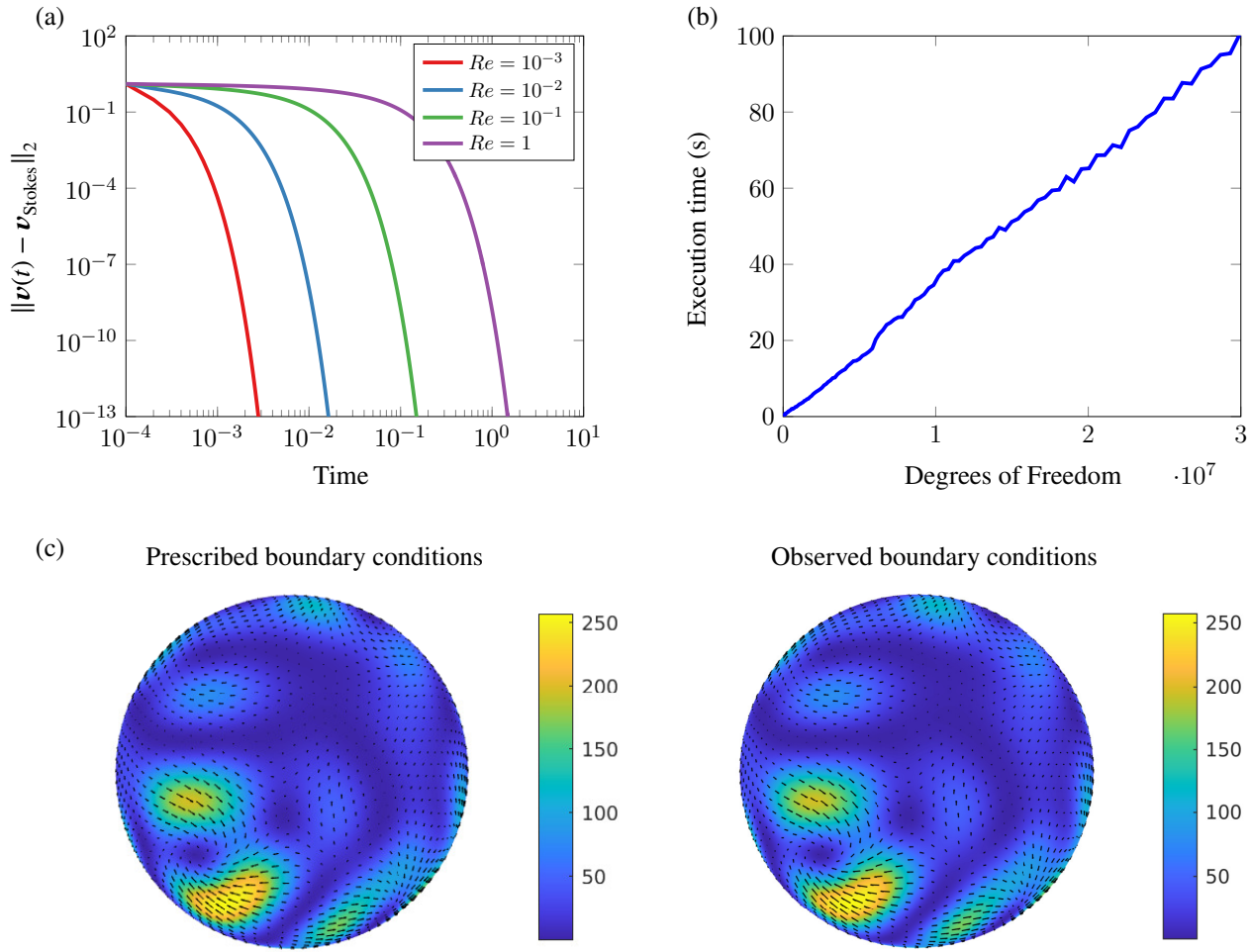
\begin{figure}[ht]
\centering
\begin{overpic}[width=0.9\textwidth]{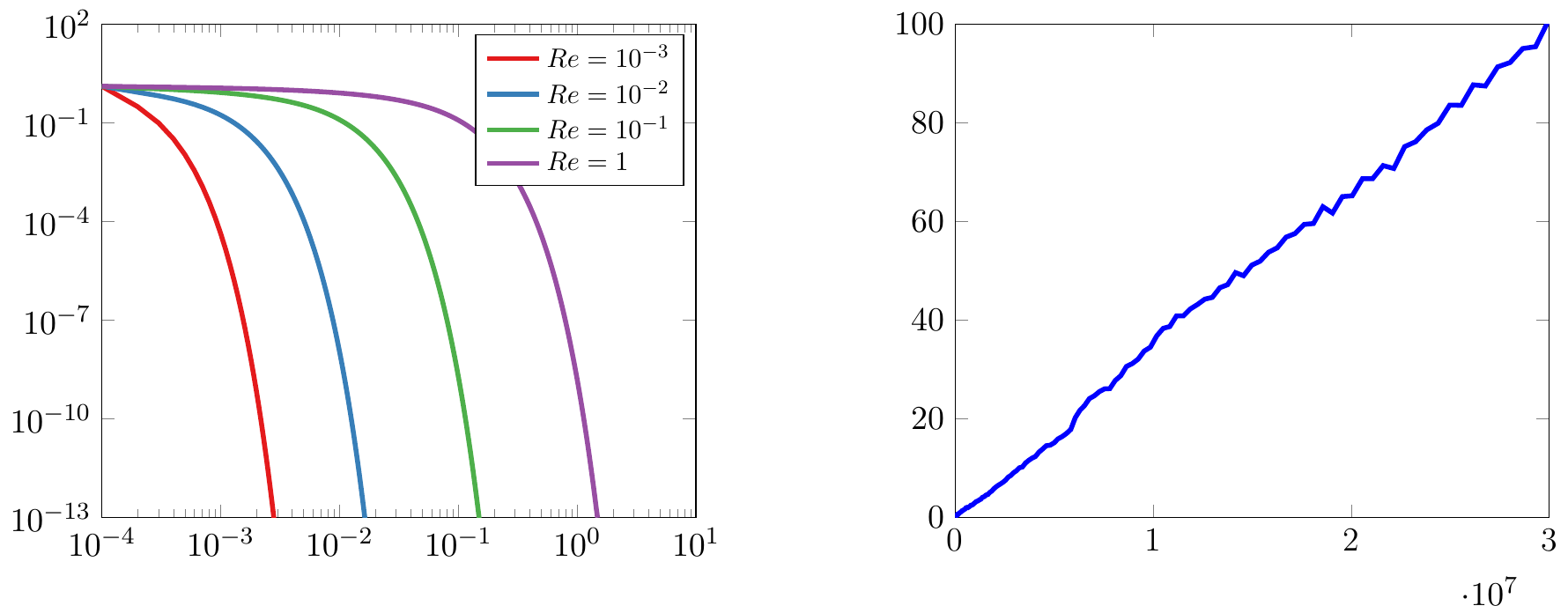}
\put(-3,14){\rotatebox{90}{$\|\v(t)-\v_{\textup{Stokes}}\|_2$}}
\put(22,0){Time}
\put(53,13){\rotatebox{90}{Execution time (s)}}
\put(71,0){Degrees of Freedom}
\put(-1,38){(a)}
\put(53,38){(b)}
\put(-1,-6){(c)}
\end{overpic}\\
\vspace{1.5cm}
\hspace{0.5cm}
\begin{overpic}[width=0.4\textwidth,trim={130 60 50 70},clip]{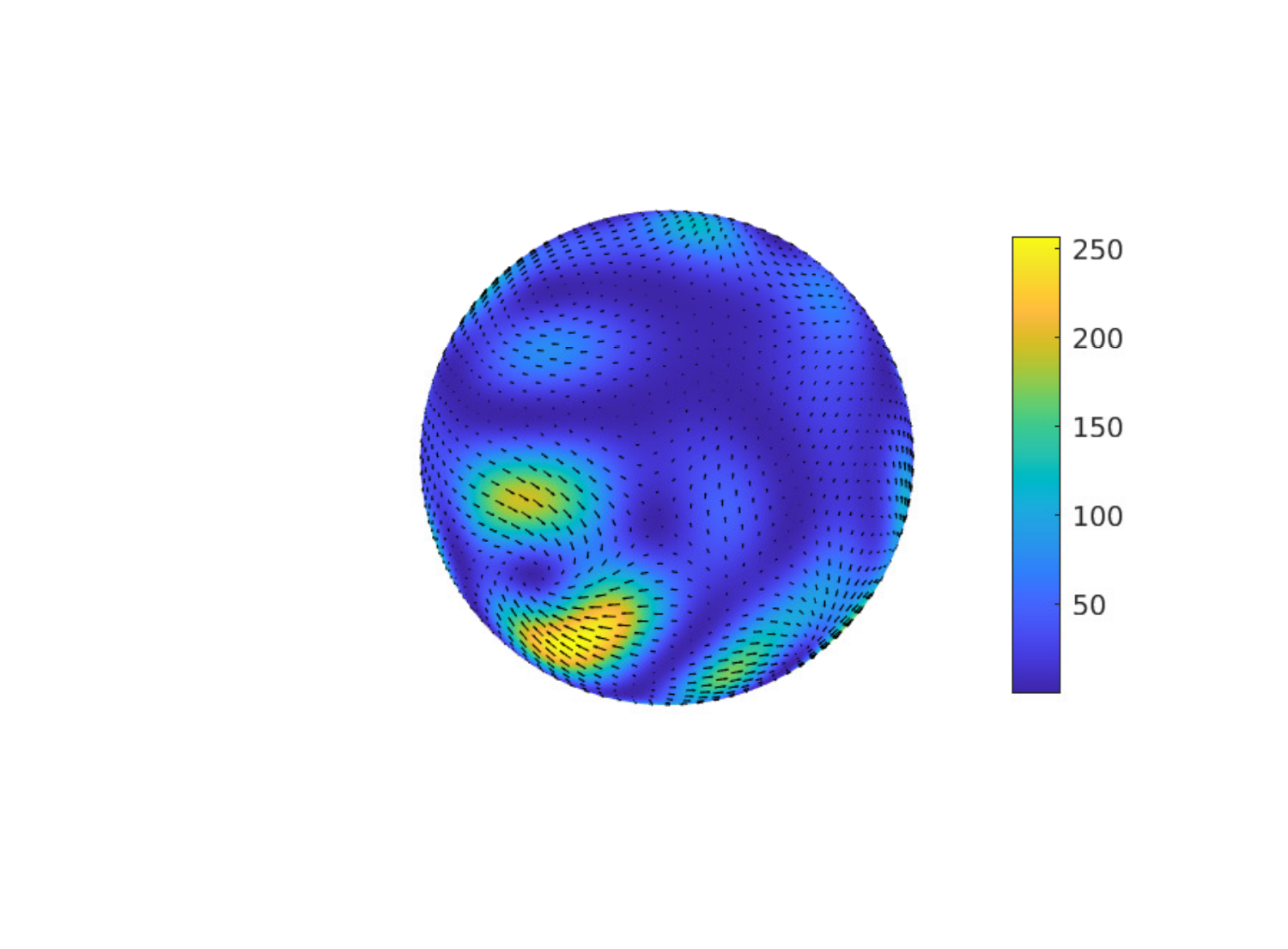}
\put(10,80){Prescribed boundary conditions}
\end{overpic}
\hspace{1cm}
\begin{overpic}[width=0.4\textwidth,trim={130 60 50 70},clip]{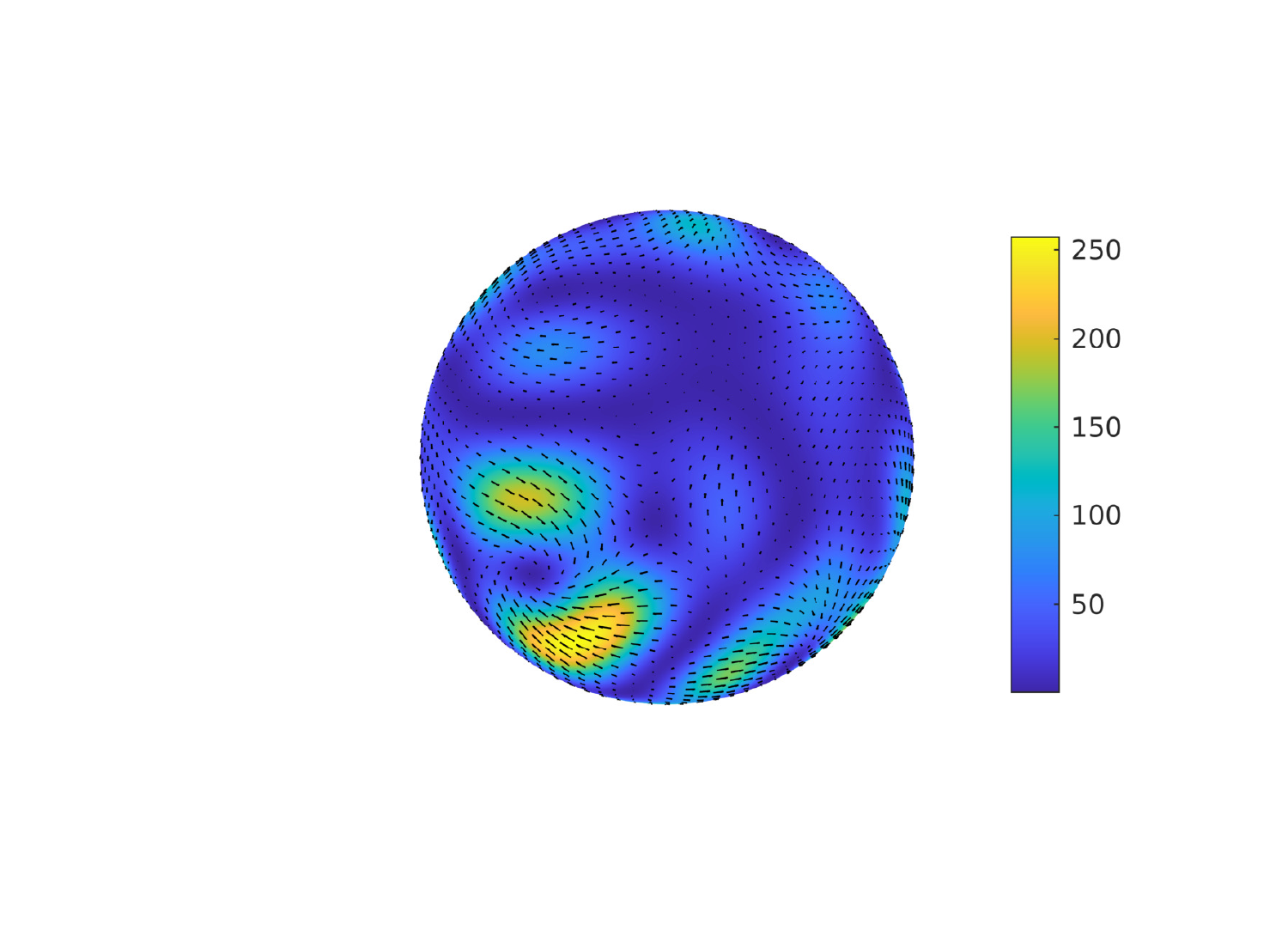}
\put(10,80){Observed boundary conditions}
\end{overpic}
\caption{Convergence of the velocity field $\v(t)$ to the analytical solution $\v_{\textup{Stokes}}=r\sin(\theta)\hat{\lambda}$ with respect to time for different Reynolds number (a). Computational timings of the algorithm for one time step with respect to the number of degrees of freedom (b). Test of the boundary conditions for the velocity field with prescribed velocity field at the surface and observed velocity at the boundary after one thousand time-steps (c).
}
\label{fig_stokes_timings}
\end{figure}

First, we conduct a convergence test of the algorithm presented in \cref{sec_numer_method} and solve the NS equations (written in vorticity form) on the unit ball:
\begin{align*}
\nabla^2\bs{\psi} &= -\bs{\omega},\\
\frac{\partial\bs{\omega}}{\partial t}-\frac{1}{Re}\nabla^2\bs{\omega} &= -\nabla\times(\bs{\omega}\times \v),
\end{align*}
for low Reynolds numbers, and boundary conditions for the velocity field given as
\[f=0,\qquad g=\cos\theta,\]
where 
\begin{equation} \label{eq_BC_NS}
\v|_{\partial\Omega}=\nabla_1 f+\Lambda_1 g = \sin(\theta)\hat{\lambda},
\end{equation}
and $\hat{\lambda}$ denotes the unit vector in the azimuthal direction. We use a smooth random initial vorticity vector field sampled from a Gaussian process with squared-exponential covariance kernel~\cite{filip2019smooth,rasmussen2006gaussian}, a spatial discretization of $n=100$, corresponding to $1.5\times 10^6$, and a time-step of $\Delta_t = 10^{-4}$. As $Re\to 0$, the velocity should converge to the velocity solution to the Stokes' equations with similar boundary conditions, given by $\v_{\text{Stokes}}=r\sin(\theta)\hat{\lambda}$. We perform 4 transient simulations of the NS equations with respective Reynolds number $Re \in\{1,10^{-1},10^{-2},10^{-3}\}$ and report the error in the $L^2$-norm between the velocity $v(t)$ and $\v_{\text{Stokes}}$ in \cref{fig_stokes_timings}(a). We observe that the error decays to nearly machine precision as the simulation time $t$ increases and that the velocity obtained converges to the solution to the Stokes's equation when $Re\to 0$ as expected. This behaviour shows the spectral accuracy of the solver resulting from the choice of the spatial discretization.

Next, we perform several NS simulations with Reynolds number $Re=70$, time-step $\Delta_t=10^{-3}$, and boundary conditions given by \cref{eq_BC_NS}. We vary the spatial discretization from $N=11$ to $n=271$, corresponding to a maximum number of degrees of freedom of $n=3\times 10^7$, and record the execution time\footnote{Timings were performed on an Intel Xeon CPU E5-2667 v2 @ 3.30GHz using MATLAB R2019b without explicit parallelization.} taken by the algorithm over the first 15 time-steps. \Cref{fig_stokes_timings}(b) displays the average computational time for one time-step as we increase the number of degrees of freedom and illustrate the optimal complexity of the method. Moreover, the NS equations can be solved on a ball with one million degrees of freedom in a dozen seconds on a single CPU core. Finally, we highlight that the numerical algorithm presented in \cref{sec_numer_method} fully decouple the NS equations into poloidal and toroidal equations, and all the subsequent operations (such as solving Helmholtz's equations, changing basis, and vector calculus) can be parallelized along one spatial direction, leading to nearly optimal scaling when using multiple CPU cores.

Finally, we test the imposition of the boundary conditions on the velocity field by performing a simulation at Reynolds number $Re=1$, with time-step $\Delta_t=10^{-3}$, a spatial discretization of $n=101$, and smooth random surface potentials $f$ and $g$ (sampled from a Gaussian process). The left panel of \cref{fig_stokes_timings}(c) displays the prescribed velocity field: $\v_{\text{prescribed}}=\nabla_1 f+\Lambda_1 g$ and its magnitude at the surface of the sphere. In contrast, the right panel consists of the simulated velocity at the boundary after one thousand time-step. We observe that $L^2$-norm of the difference between the two fields at the surface of the sphere is bounded by $10^{-8}$.

\subsection{Active fluids simulation}

We perform a simulation of the active fluid model in Eq.~(\ref{eq_active_fluids}) with dimensionless parameters $\Gamma_0=1,\, \Gamma_2=-8.13\times 10^{-3},\, \Gamma_4=1.65\times 10^{-5}$. This choice of parameters represents active energy injection into the fluid at the scale corresponding to vortices with characteristic size $\Gamma=\pi (-2\Gamma_4/\Gamma_2)^{1/2}=0.2$ and growth timescale $\tau=[\Gamma_2/(2\Gamma_4)(\Gamma_0-\Gamma_2^2/(4\Gamma_4))]^{-1}=2.76$ confined to a narrow active bandwidth $\kappa \Lambda=0.12$, where $\kappa=[-\Gamma_2/\Gamma_4-2\sqrt{\Gamma_0/\Gamma_4}]^{1/2}$; the three characteristic scales $(\Lambda,\tau,\kappa)$ fully determine the resulting dynamics of the GNS model~\cite{slomka2017spontaneous}. As the boundary conditions for the velocity field on the unit ball, we use the no-slip boundary conditions $\v|_{\partial\Omega}=0$, which corresponds to setting the velocity surface potentials to zero $f=g=0$~(Section~\ref{sec_Integral_T_omega}). To close the 6-th order Generalized Navier--Stokes  Eq.~(\ref{eq_active_fluids}), we additionally set the Laplacian and bi-Laplacian of the toroidal and poloidal scalars of the vorticity to zero on the boundary
\begin{equation}
\label{eq:ho_bc}
\nabla^2 T_{\bs{\omega}}(1,\lambda,\theta)=\nabla^4 T_{\bs{\omega}}(1,\lambda,\theta)=0,
\quad
\nabla^2 P_{\bs{\omega}}(1,\lambda,\theta)=\nabla^4 P_{\bs{\omega}}(1,\lambda,\theta)=0.
\end{equation}
This higher-order closure is motivated by similar closures on the vorticity field previously used in two-dimensional GNS systems~\cite{Slomka2017,slomka2017stokes}. In simulations, the initial vorticity field is a smooth random vector field sampled from a Gaussian process. We use a spatial discretization parameter of $n=100$, a time-step of $\Delta_t=10^{-3}$, and run the simulation until $t=75$, corresponding to $75000$ time-steps and approximately 100 hours of execution time on a single core of a Intel Xeon CPU E5-2667 v2 @ 3.30GHz. 

\begin{figure}[ht]
\centering
\begin{overpic}[width=0.4\textwidth]{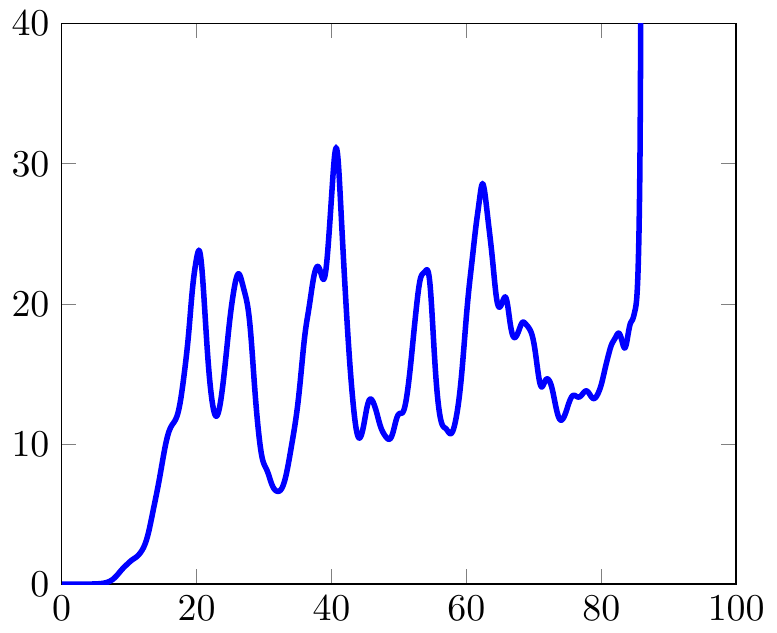}
\put(-7,27.5){\rotatebox{90}{Kinetic energy}}
\put(47,-5){Time}
\end{overpic}
\vspace{0.2cm}
\caption{Kinetic energy of the active fluid simulation with respect to the time. The energy blow-up happens at time $t \approx 86$.}
\label{fig_kinetic_energy}
\end{figure}

In a simulation, starting from a random field, the system builds up kinetic energy, defined as
\[\text{Kinetic energy} = \int_{B(0.1)}\v\cdot\v\,\textup{d}x,\]
driven by the linear term proportional to $\Gamma_2$ in Eq.~(\ref{eq_active_fluids})~(see \cref{fig_kinetic_energy}). As the kinetic energy grows, the advective nonlinearity starts to impact the dynamics by mixing different wavelengths. As a result of this nonlinear coupling, the kinetic energy growth is eventually balanced by the dissipation at large and small wavelengths due to the viscous and hyper-viscous terms proportional to $\Gamma_0$ and $\Gamma_4$. Once this balance is reached, the system settles onto a statistically stationary state characterized by vortices of well-defined size~(see\cref{fig_active_fluids} and a movie in Supplementary Material). The presence of vortices of well-defined size reflects our parameter choice to simulate a system with scale selection by settling a small value of the active bandwidth $\kappa$. At larger times, after a stationarity period, the kinetic energy grows rapidly again, and the system becomes unstable. We are unsure about the origin of this long-time instability. A possible reason is due to too small numerical discretization. Another reason for the blow-up is the boundary conditions in Eq.~(\ref{eq:ho_bc}) being not dissipative enough, allowing the system to accumulate unlimited kinetic energy.

\begin{figure}[htbp] 
\centering
\vspace{0.5cm}
\begin{overpic}[height=0.23\textwidth,trim={120 70 105 67},clip]{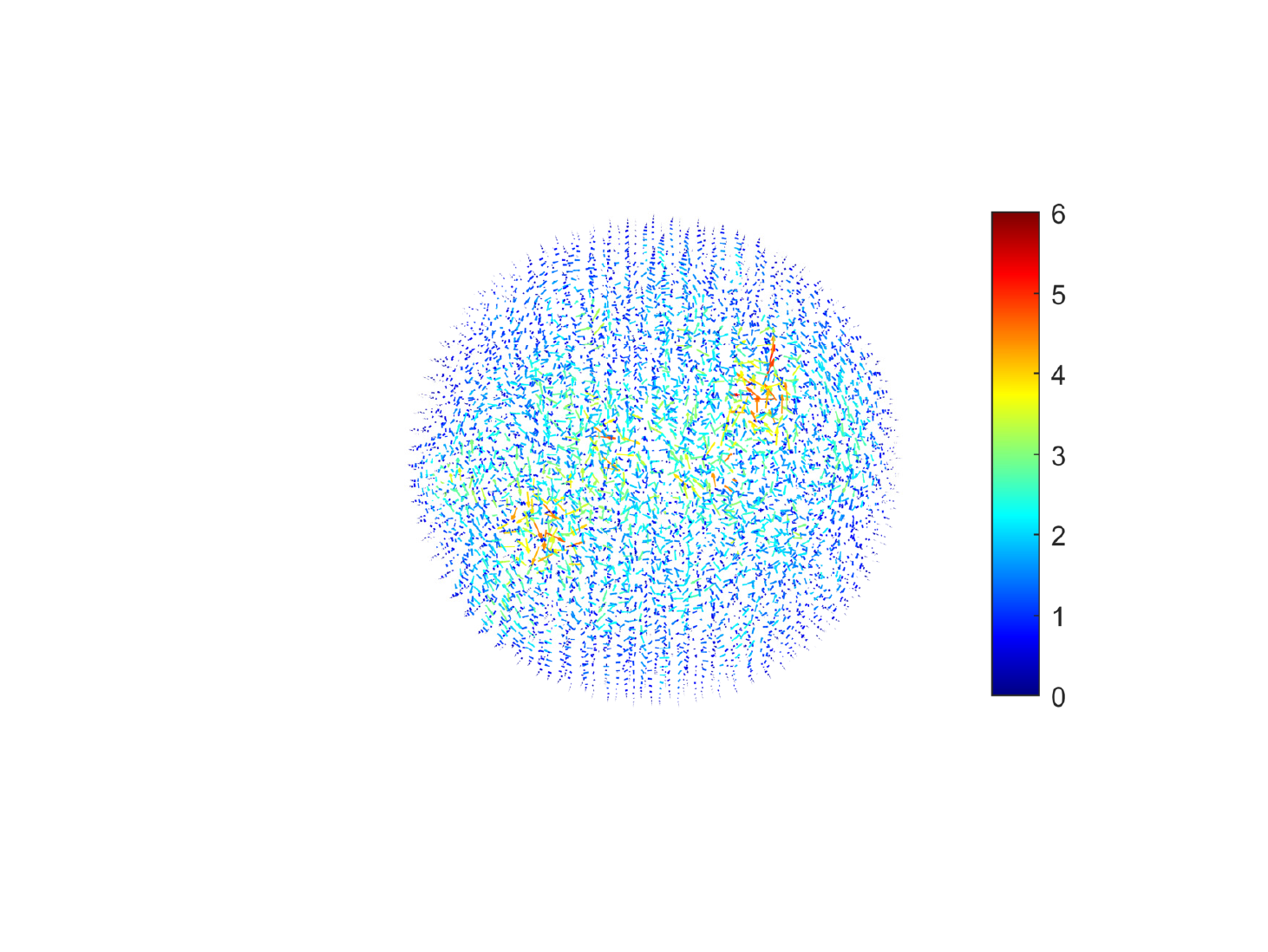}
\put(40,95){$t=50$}
\end{overpic}
\hspace{0.2cm}
\begin{overpic}[height=0.23\textwidth,trim={120 70 105 67},clip]{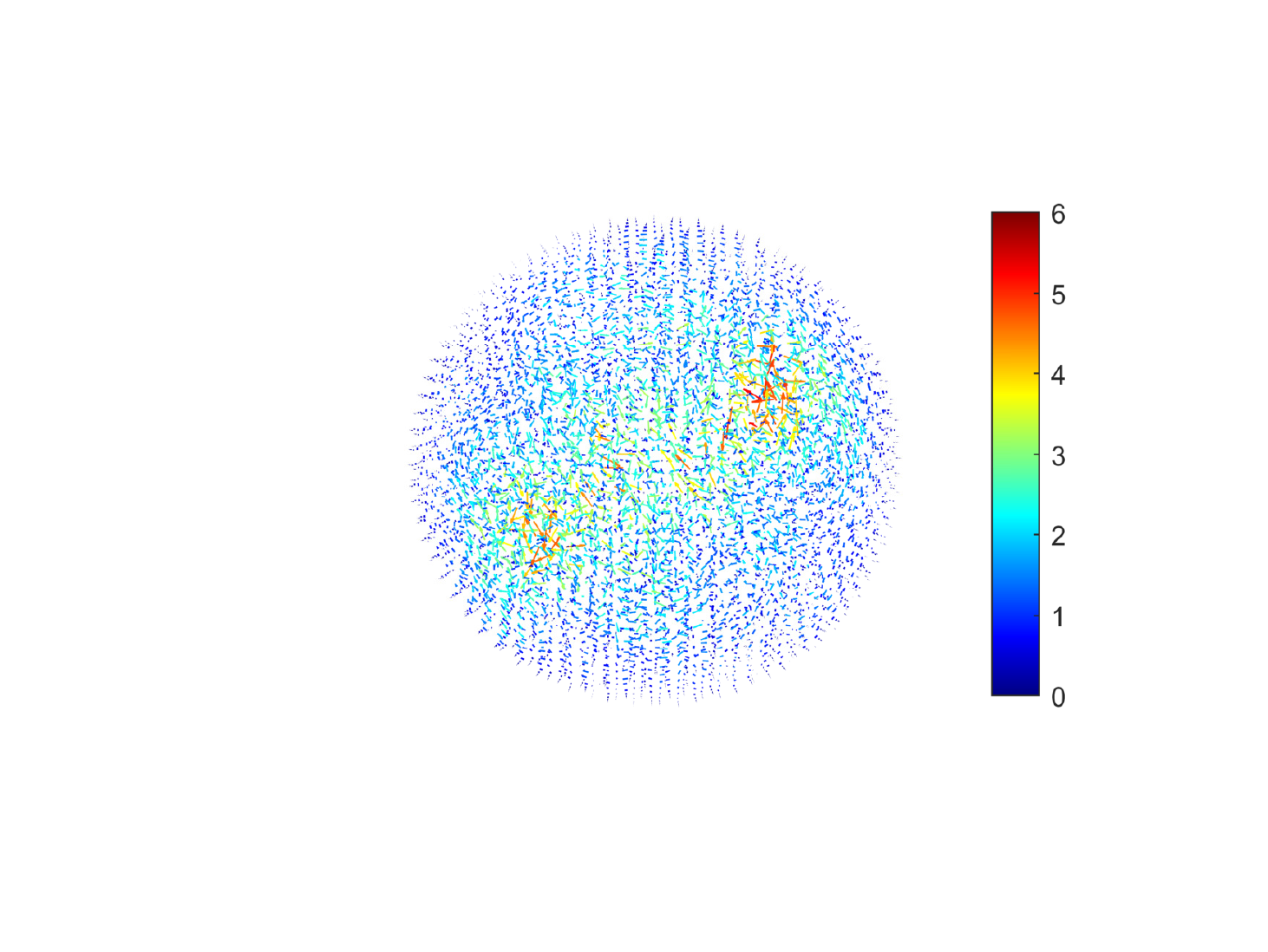}
\put(40,95){$t=51$}
\put(148,95){$t=52$}
\end{overpic}
\hspace{0.2cm}
\begin{overpic}[height=0.23\textwidth,trim={120 70 50 67},clip]{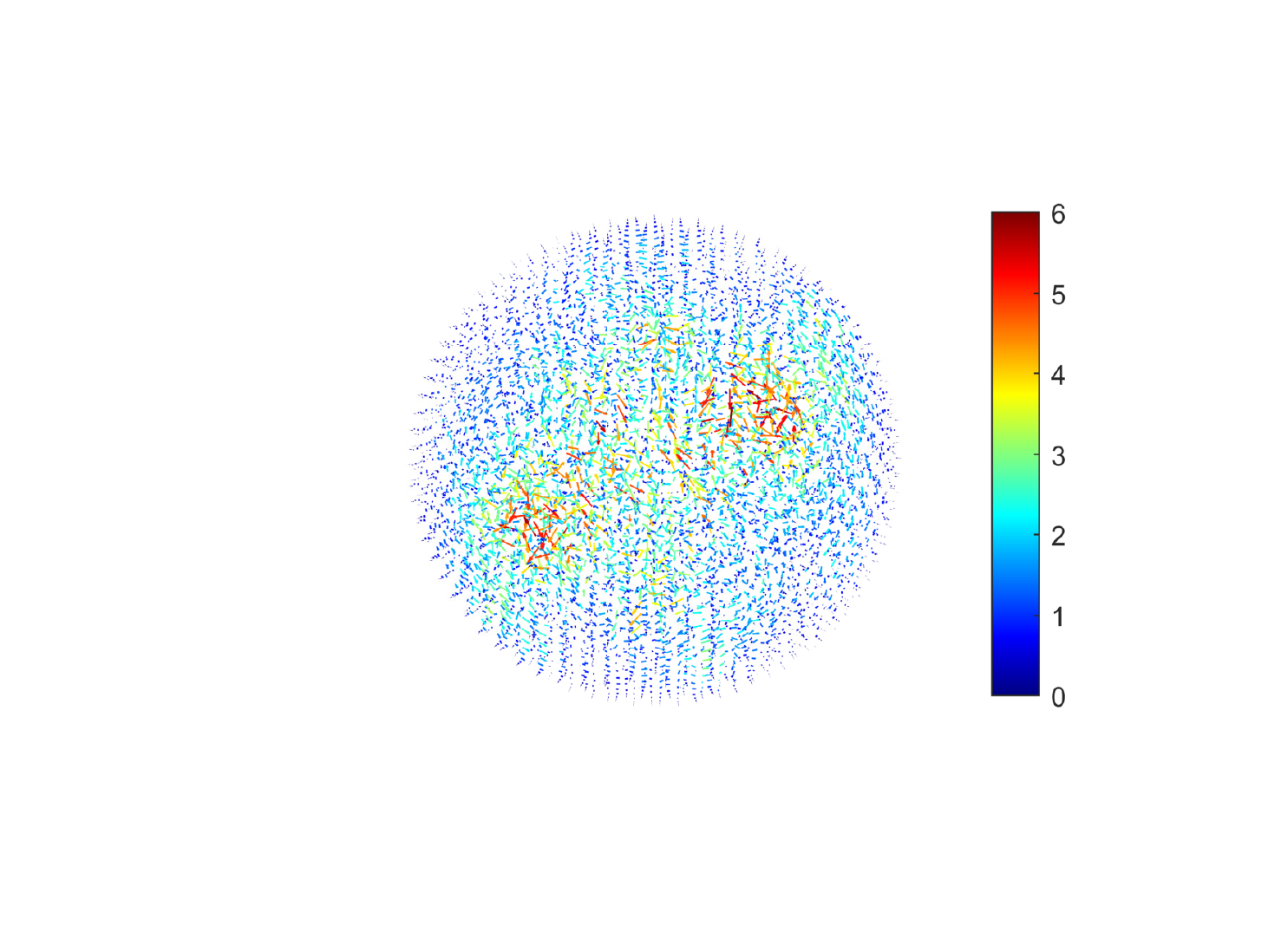}
\end{overpic}\\
\begin{overpic}[height=0.23\textwidth,trim={120 70 105 67},clip]{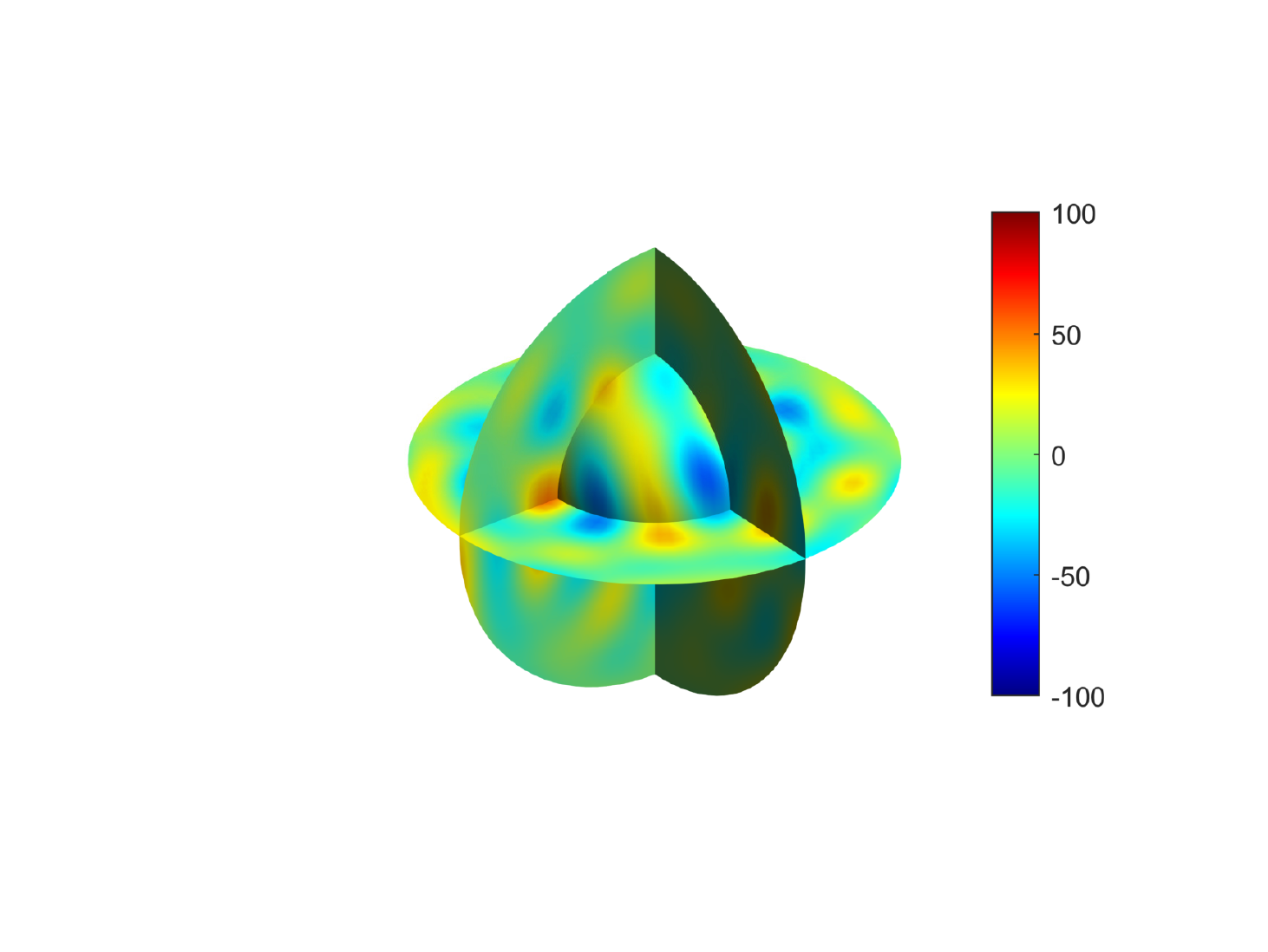}
\end{overpic}
\hspace{0.2cm}
\begin{overpic}[height=0.23\textwidth,trim={120 70 105 67},clip]{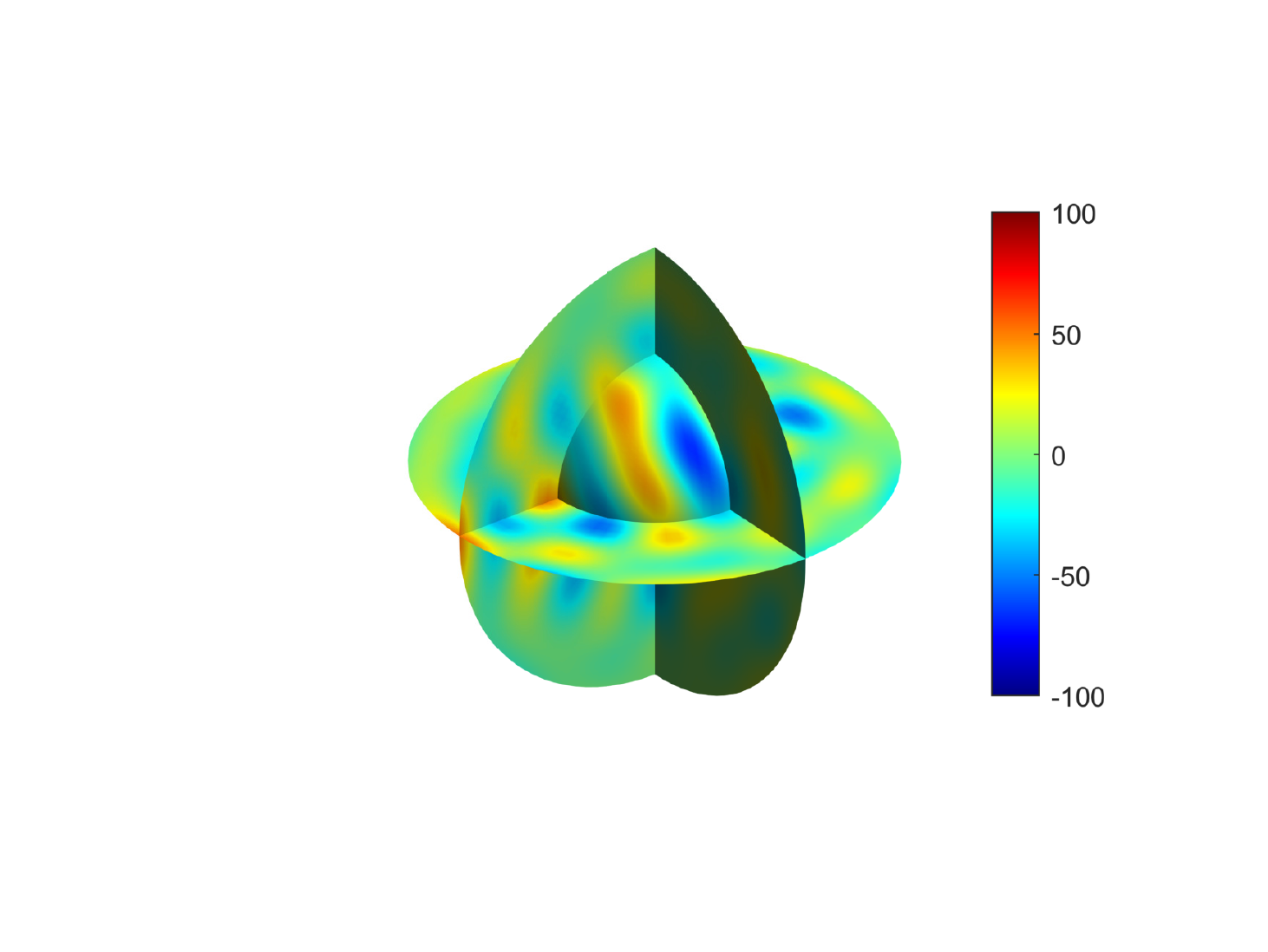}
\end{overpic}
\hspace{0.2cm}
\begin{overpic}[height=0.23\textwidth,trim={120 70 50 67},clip]{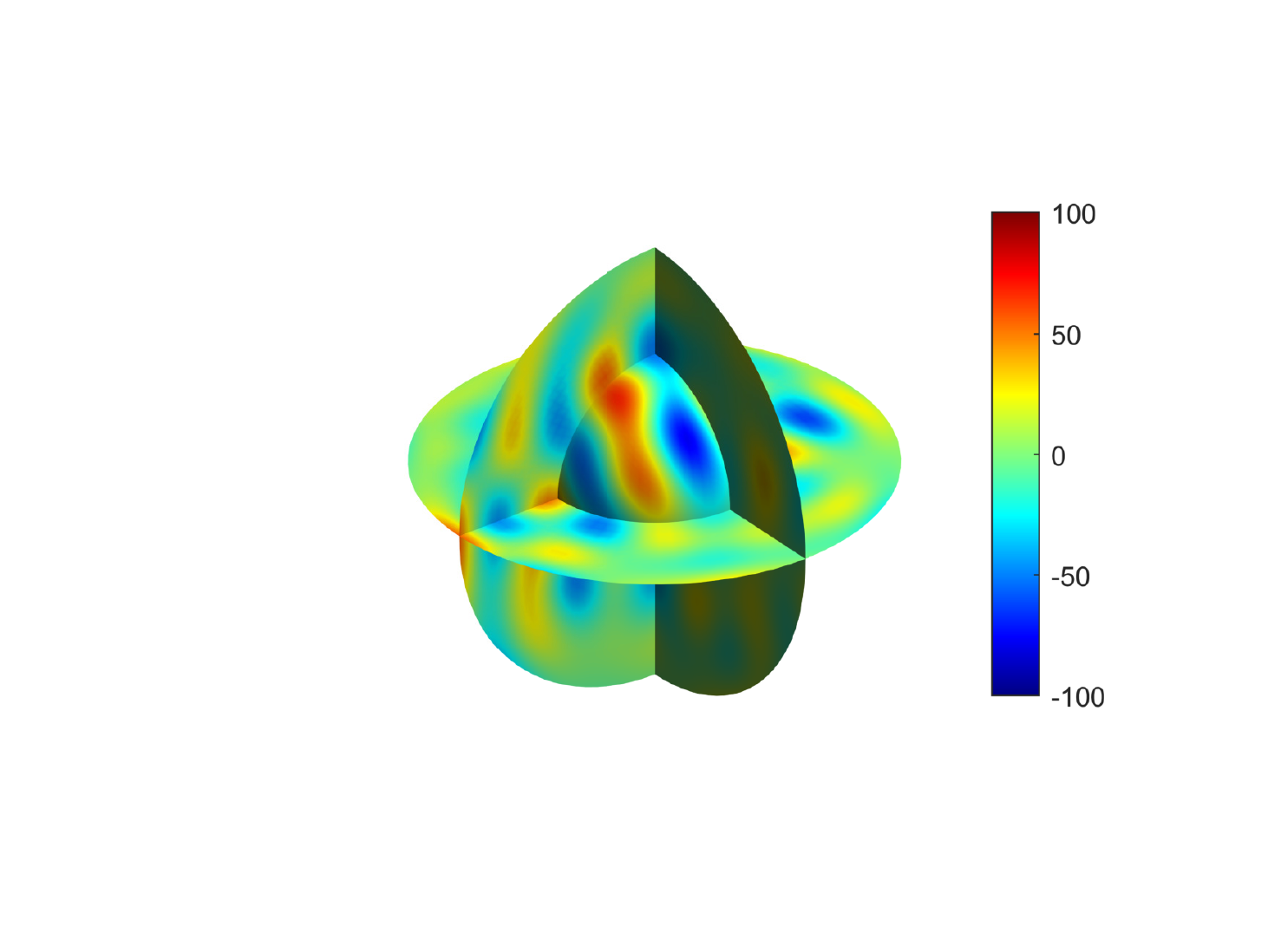}
\end{overpic}\\
\vspace{1cm}
\begin{overpic}[height=0.23\textwidth,trim={120 70 105 67},clip]{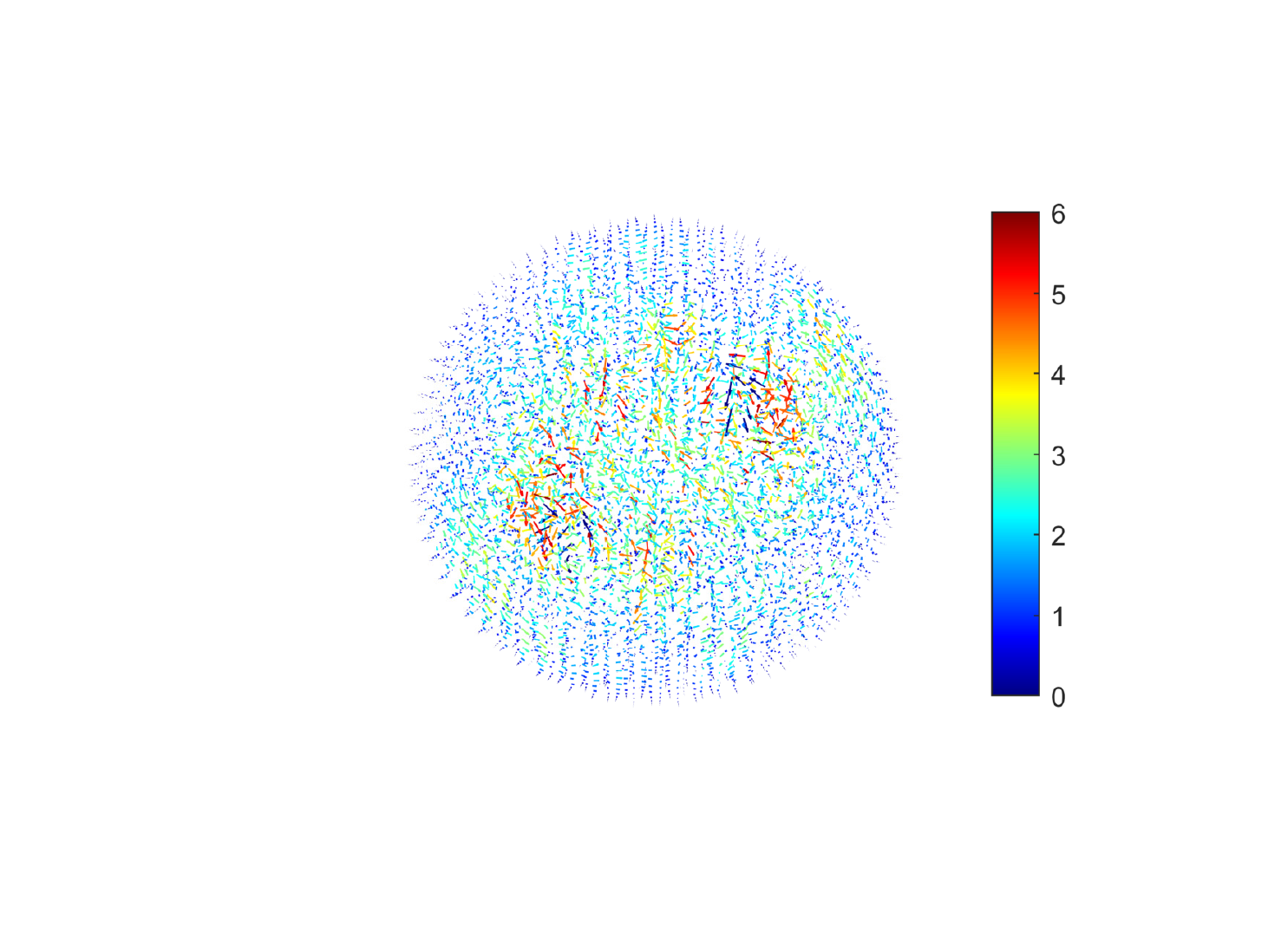}
\put(40,95){$t=53$}
\end{overpic}
\hspace{0.2cm}
\begin{overpic}[height=0.23\textwidth,trim={120 70 105 67},clip]{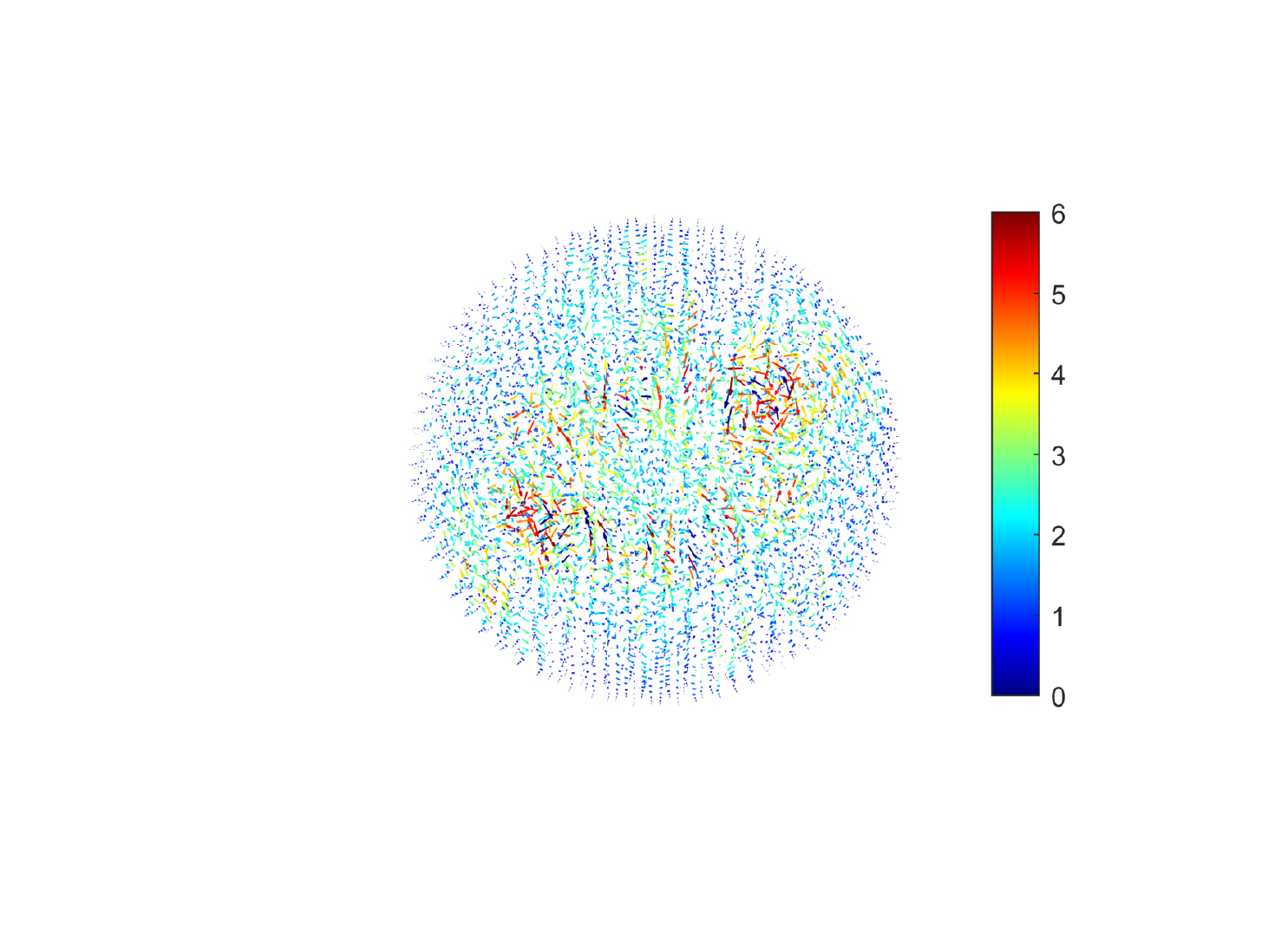}
\put(40,95){$t=54$}
\put(148,95){$t=55$}
\end{overpic}
\hspace{0.2cm}
\begin{overpic}[height=0.23\textwidth,trim={120 70 50 67},clip]{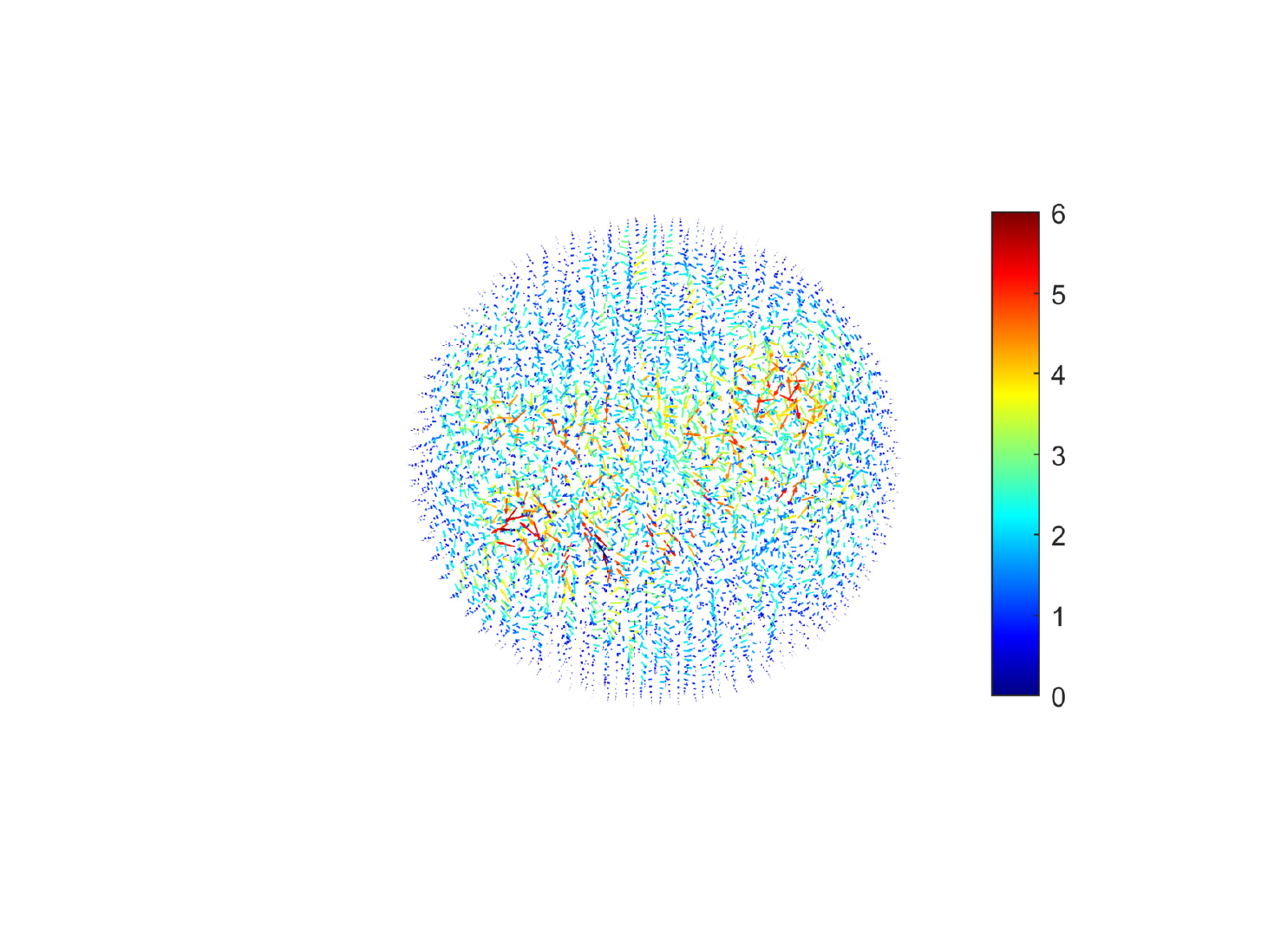}
\end{overpic}\\
\begin{overpic}[height=0.23\textwidth,trim={120 70 105 67},clip]{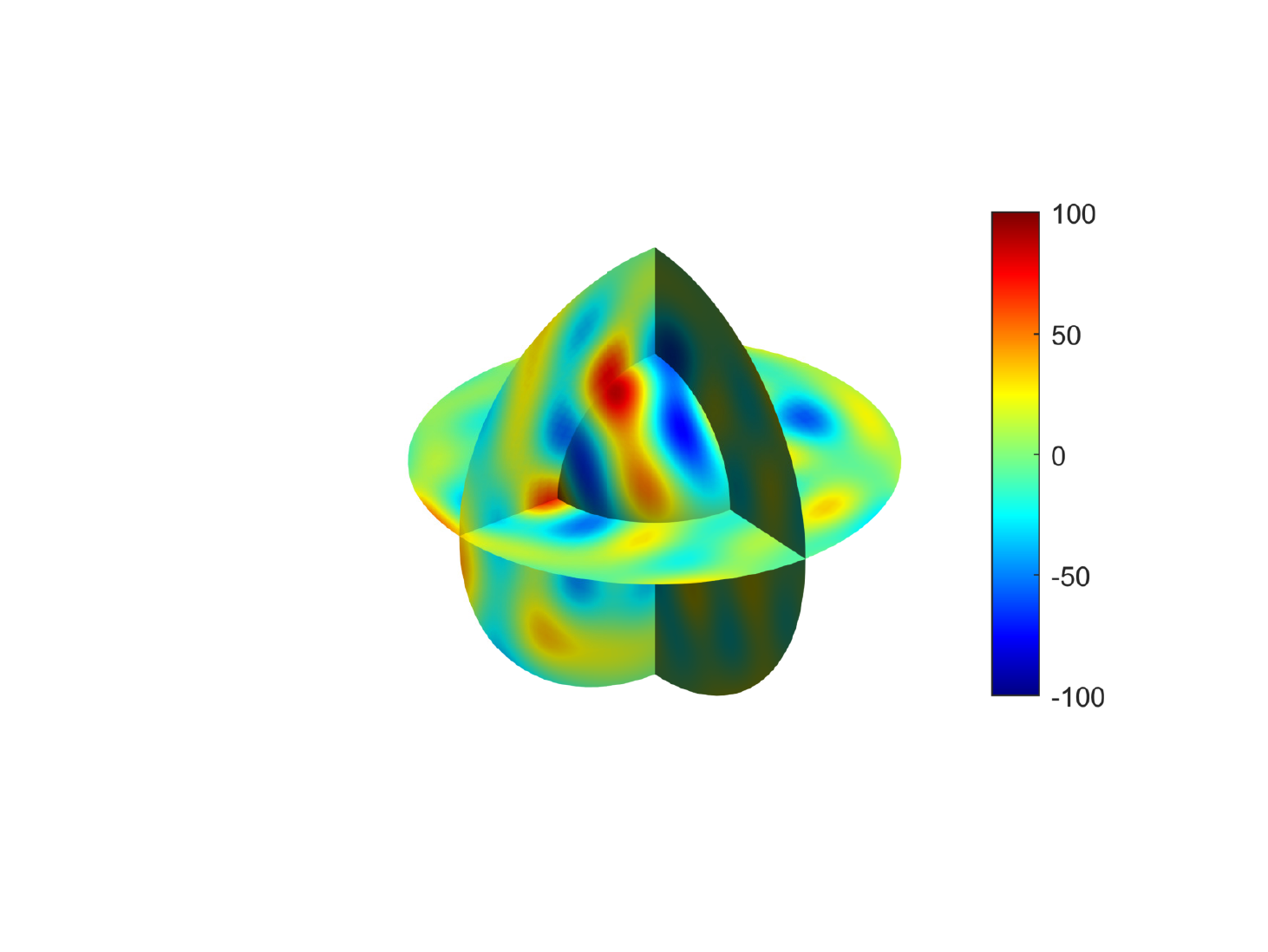}
\end{overpic}
\hspace{0.2cm}
\begin{overpic}[height=0.23\textwidth,trim={120 70 105 67},clip]{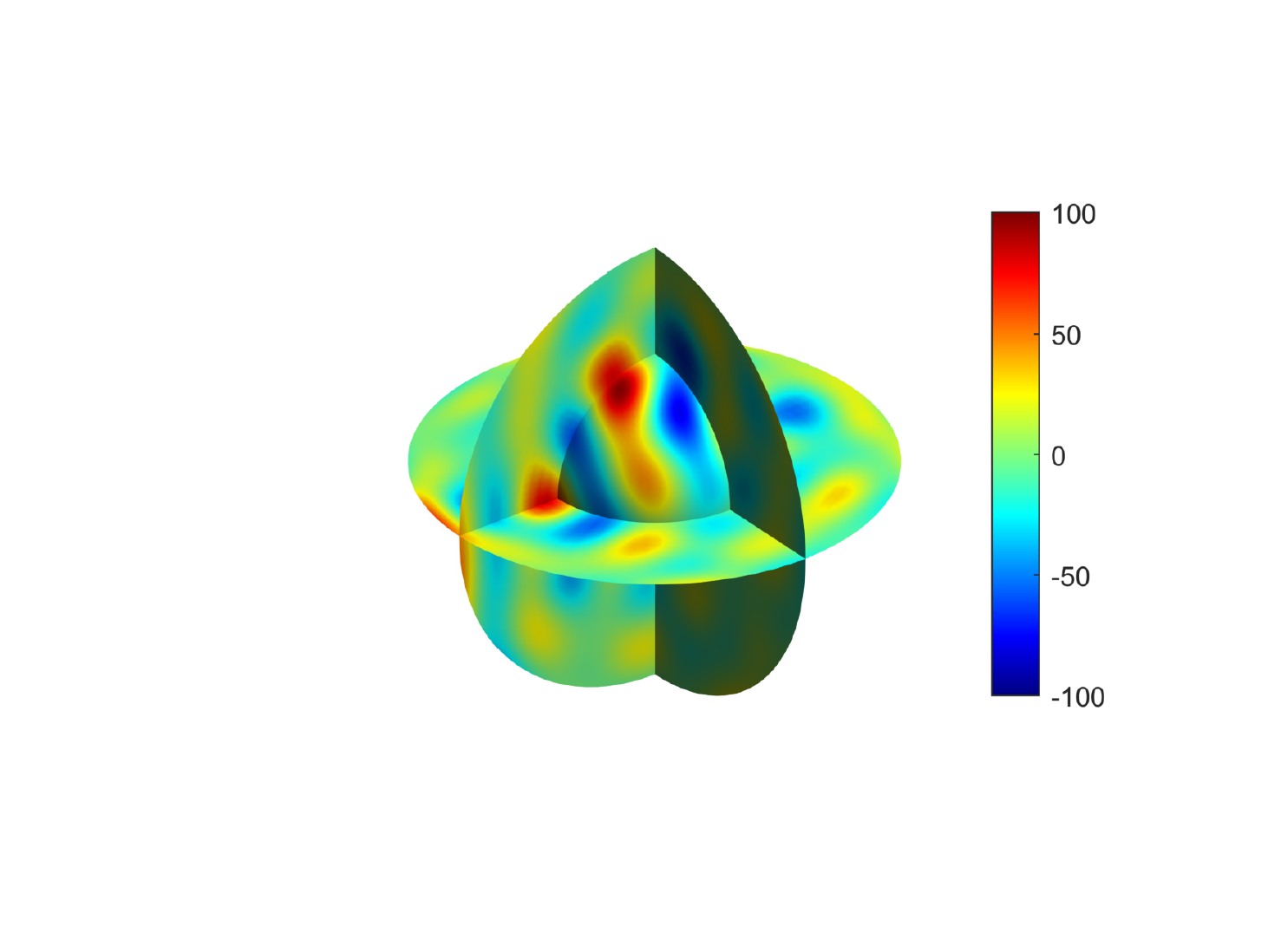}
\end{overpic}
\hspace{0.2cm} 
\begin{overpic}[height=0.23\textwidth,trim={120 70 50 67},clip]{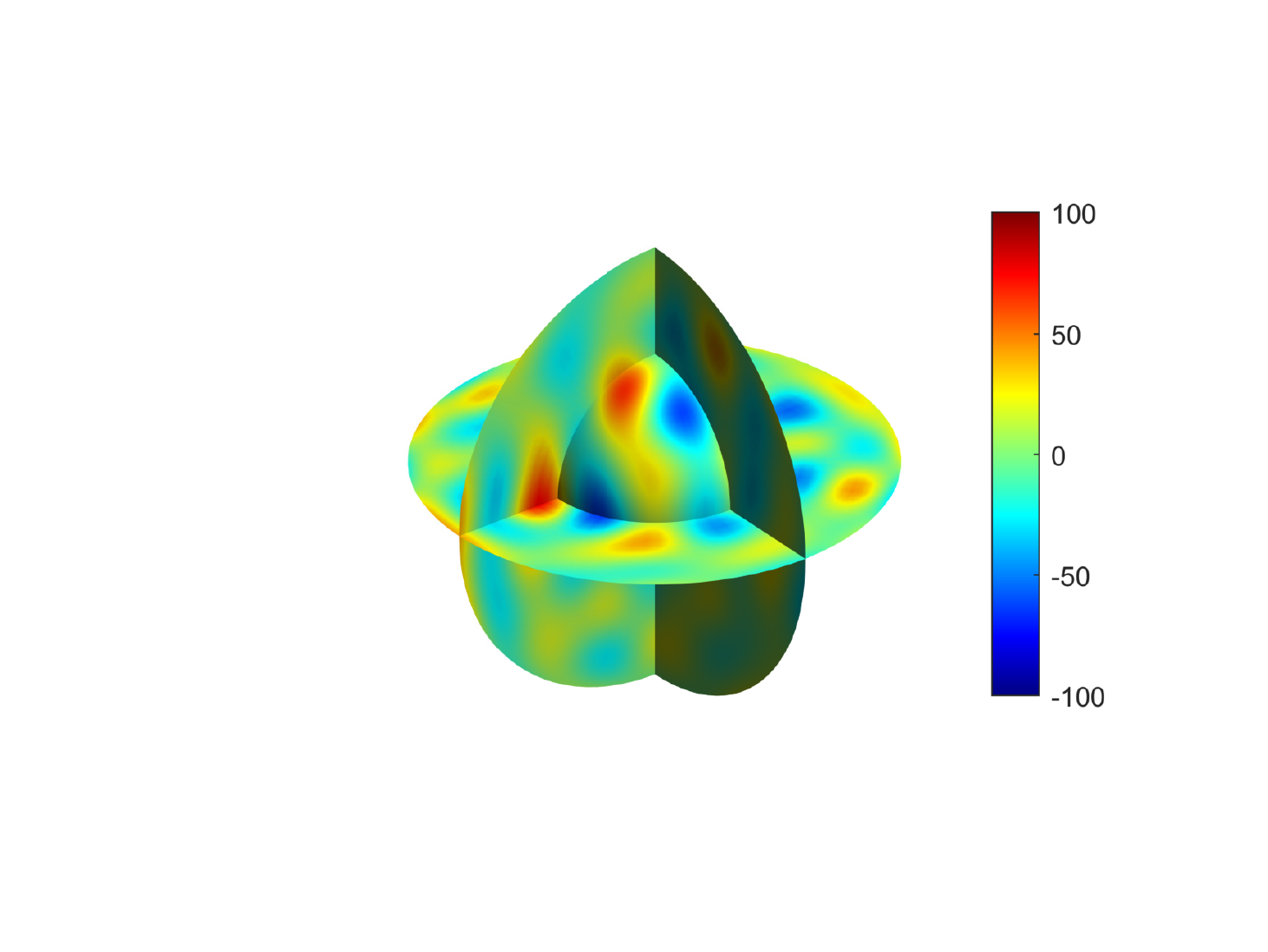}
\end{overpic}

\caption{Snapshots of the velocity field and $z$ component of the vorticity field during the active fluids simulation at times $t=50,\ldots,55$.}
\label{fig_active_fluids}
\end{figure}

\section{Conclusions} \label{sec_concl}

By exploiting the poloidal-toroidal decomposition of divergence-free vector fields, the Chebyshev--Fourier--Fourier and Chebyshev-Spherical harmonics bases, we developed an algorithm for solving the incompressible generalized NS equations on the ball with spectral accuracy and optimal complexity per time-step. Numerical experiments illustrated the linear scaling of the execution time with respect to the number of degrees of freedom and applied it to simulate an active fluids problem.

Our algorithm makes it straightforward to couple any surface dynamics with the interior dynamics once the user specifies separate evolution equations for the surface potentials of the compressible and rotational components of the velocity field on the boundary. Such active surface-driven flows could be important, for example, in studies of bulk flows driven by active stresses confined to a boundary~\cite{Mickelin2018,Shankar2017,Supekar2020,Rank2021}.

The instability of the active fluid simulation and blow-up of the corresponding kinetic energy at large time can be investigated in the future by considering alternative higher-order boundary conditions to introduce more dissipation in the system. Another extension of this work is to parallelize the current implementation of the numerical method presented in this paper in order to reach the large spatial discretization required to simulate Navier--Stokes turbulences at large Reynolds number and study the stability of the technique in this regime.

\section*{Acknowledgements}

We thank J\"orn Dunkel for discussions and Keaton Burns for making us aware of performance benchmarks~\cite{marti2014full,matsui2016performance}. This work was supported by the EPSRC Centre for Doctoral Training in Industrially Focused Mathematical Modelling (EP/L015803/1) in collaboration with Simula Research Laboratory. The third author was supported by the National Science Foundation grant DMS-1818757, DMS-1952757, and DMS-2045646, as well as the FACE Foundation.

\bibliography{biblio}%

\end{document}